\documentclass{article}%
\usepackage{latexsym}
\usepackage{amsmath}
\usepackage{amsfonts}
\usepackage{amssymb}
\usepackage{graphicx}%
\providecommand{\U}[1]{\protect\rule{.1in}{.1in}}
\newtheorem{theorem}{Theorem}
\newtheorem{acknowledgement}[theorem]{Acknowledgement}

\newtheorem{conjecture}[theorem]{Conjecture}
\newtheorem{corollary}[theorem]{Corollary}

\newtheorem{definition}[theorem]{Definition}

\newtheorem{lemma}[theorem]{Lemma}

\newtheorem{problem}[theorem]{Problem}
\newtheorem{proposition}[theorem]{Proposition}
\newtheorem{remark}[theorem]{Remark}

\begin{document}

\author{Andrey Todorov\thanks{This work was done during the stay of the author at
Max-Plank Institute f\"{u}r Mathematik-Bonn. The author wants to thank
Max-Plank Institute f\"{u}r Mathematik-Bonn for the financial support and
creation of extremely nice conditions for work. }\\University of California\\Department of Mathematics,\\Santa Cruz, CA 95064\\Bulgarian Academy of Sciences\\Institute of Mathematics\\Sofia, Bulgaria}
\title{The Analogue of Dedekind Eta Functions for Calabi-Yau Manifilolds II. \\(Algebraic, Analytic Discriminants and the Analogue of Baily-Borel
Compactification of the Moduli Space of CY Manifolds.) }
\maketitle
\date{}

\begin{abstract}
In this paper we construct the analogue of Dedekind $\eta-$function on the
moduli space of polarized CY manifolds. We prove that the $\mathbf{L}^{2}$
norm of $\eta(\tau)$ is the regularized determinants of the Laplacians of the
CY metric on $(0,1)$ forms.

We construct the analogue of the Baily-Borel Compactification of the moduli
space of polarized CY and prove that it has the same properties as the
Baily-Borel compactification of the locally symmetric Hermitian spaces. We
proved that the compactification constructed in the paper is the minimal.

\end{abstract}
\tableofcontents

\section{Introduction}

This is the second part of the paper \cite{BT}. In this paper we will study
the algebraic and analytic discriminants of Calabi-Yau manifolds. We will
generalize the notion of algebraic and analytic discriminants of elliptic
curves to the case of CY manifolds. Next we will review the notion of
algebraic and analytic discriminants for elliptic curves and their relation.

\textbf{1. Algebraic Discriminant on the Moduli Space of Marked Elliptic
Curves.}

The algebraic discriminant of an elliptic curves is defined on the moduli
space of elliptic curves which is the quotient of the Teichm\"{u}ller space by
the mapping class group. The Teichm\"{u}ller space is the upper half plane
$\mathfrak{h}:=\left\{  \tau\in\mathbb{C}|\operatorname{Im}\tau>0\right\}  $
and the mapping class group of the elliptic curve is $\mathbb{SL}_{2}\left(
\mathbb{Z}\right)  .$ Thus the moduli space of elliptic curves is
$\mathbb{PSL}_{2}\left(  \mathbb{Z}\right)  \left\backslash \mathfrak{h}%
\right.  .$ The elliptic curves corresponding to lattices spanned by
$(1,\sqrt{-1})$ and $(1,\rho),$ where $\rho^{3}=1$ and $\rho\neq1$ have
automorphisms of order four and six. So the moduli space $\mathbb{PSL}%
_{2}\left(  \mathbb{Z}\right)  \left\backslash \mathfrak{h}\right.  $ is a
stack. If we consider the moduli space of elliptic curves with a fixed basis
in $H_{1}\left(  E,\mathbb{Z}/2\mathbb{Z}\right)  ,$ then it is isomorphic to
$\Gamma(2)\left\backslash \mathfrak{h}\right.  ,$ where
\[
\Gamma(2):=\left\{  A=\left(
\begin{array}
[c]{cc}%
a & b\\
c & d
\end{array}
\right)  \left\vert \det A=1;\text{ }a,b,c,d\in\mathbb{Z},\right.
A\equiv\left(
\begin{array}
[c]{cc}%
1 & 0\\
0 & 1
\end{array}
\right)  \operatorname{mod}2\right\}  .
\]
The group $\Gamma(2)$ acts on $\mathfrak{h}$ without fixed points. Thus over
$\Gamma(2)\left\backslash \mathfrak{h}\right.  $ we have a universal family of
elliptic curves with a fixed basis in $H_{1}\left(  E,\mathbb{Z}%
/2\mathbb{Z}\right)  .$ This universal family can be represented by
\begin{equation}
y^{2}=x\left(  x-1\right)  \left(  x-\lambda\right)  . \label{1}%
\end{equation}
So $\Gamma(2)\left\backslash \mathfrak{h}\right.  =\mathbb{CP}^{1}-\left\{
0,1,\infty\right\}  .$ The algebraic discriminant $\Delta_{alg}(\tau)$ of the
family $\left(  \ref{1}\right)  $ is given by
\begin{equation}
\Delta_{alg}(\lambda)=\left(  \left(  1-\lambda\right)  \lambda\right)  ^{2}.
\label{2}%
\end{equation}

\textbf{2. Analytic Discriminant.}

The analytic discriminant of the elliptic curve $E_{\tau}:=\mathbb{C}\left/
\left\{  m+n\tau\right\}  \right.  $ is just the regularized determinant of
the Laplacian $\Delta_{(0,1)}(\tau)$ of the flat metric acting on $(0,1)$
forms. It is defined as follows: Let $0<\lambda_{1}\leq...\leq\lambda_{k}%
\leq...$ be the spectrum of $\Delta_{(0,1)}(\tau).$ Let $\zeta_{\Delta
_{(0,1)}}(\tau,s):=%
{\displaystyle\sum\limits_{i=1}^{\infty}}
\frac{1}{\lambda_{i}^{s}}.$ It is a well known fact that $\zeta_{\Delta
_{(0,1)}}(\tau,s)$ is a meromorphic function on $\mathbb{C}$ well defined at
$0.$ Then
\begin{equation}
\Delta_{(0,1)}(\tau):=\exp\left(  \left.  \left(  -\frac{d}{ds}\zeta
_{\Delta_{(0,1)}}(\tau,s)\right)  \right\vert _{s=0}\right)  . \label{rd}%
\end{equation}
It is easy to see that the regularized determinant $\Delta_{(0,1)}(\tau)$ is a
function on $\mathfrak{M}_{2}.$ The Kronecker limit formula gives an explicit
expression of the regularized determinant:

\begin{theorem}
\label{kron}Let $\det\Delta_{\tau}$ be the regularized determinant of the flat
metric on the elliptic curve $E_{\tau}:=\mathbb{C}/\left\{  n+m\tau\right\}
.$ Then $\det\Delta_{\tau}=\operatorname{Im}\tau\left\vert \eta(\tau
)\right\vert ^{2},$ $\eta(\tau)$ is defined as: $\eta(\tau)=q^{1/24}%
{\displaystyle\prod\limits_{n=1}^{\infty}}
\left(  1-q^{n}\right)  ,$ where $q=e^{2\pi i\tau}.$
\end{theorem}

\textbf{3. Relations between the Analytic and Algebraic Discriminants of the
Elliptic Curves. }

On the moduli space $\mathfrak{M}_{2}:=\Gamma(2)\backslash\mathfrak{h}$ of
marked $H_{1}(E,\mathbb{Z}/2\mathbb{Z})$ elliptic curves, the algebraic
discriminant is a section of the line bundle $\mathcal{L}$ on $\mathfrak{M}%
_{2}:=\Gamma(2)\backslash\mathfrak{h},$ associated with the principle bundle%
\[
\mathbb{U}(1)\rightarrow\mathbb{SL}_{2}\left(  \mathbb{R}\right)
\rightarrow\mathbb{SL}_{2}\left(  \mathbb{R}\right)  \left/  \mathbb{U}%
(1)\right.  =\mathfrak{h}.
\]
It is well known that $\mathfrak{M}_{2}$ is isomorphic to $\mathbb{CP}%
^{1}-\left\{  0,1,\infty\right\}  .$ The universal cover of $\mathbb{CP}%
^{1}-\left\{  0,1,\infty\right\}  $ is the upper half plane $\mathfrak{h}$.
The fundamental group $\pi_{1}\left(  \mathbb{CP}^{1}-\left\{  0,1,\infty
\right\}  \right)  $ of $\mathbb{CP}^{1}-\left\{  0,1,\infty\right\}  $ is the
free group with two generators isomorphic to $\Gamma(2)$. There is a universal
family of marked $H_{1}(E,\mathbb{Z}/2\mathbb{Z})$ elliptic curves
\begin{equation}
\pi:\mathcal{E}\rightarrow\mathbb{CP}^{1}-\left\{  0,1,\infty\right\}  ,
\label{10}%
\end{equation}
given by $y^{2}=x\left(  x-1\right)  \left(  x-\lambda\right)  .$ The line
bundle $\mathcal{L}$ on $\mathfrak{M}_{2}:=\Gamma(2)\backslash\mathfrak{h}$ is
isomorphic to the dual of $\pi_{\ast}\omega_{\mathcal{E}/\mathfrak{M}_{2}},$
where $\omega_{\mathcal{E}/\mathfrak{M}_{2}}$ is the relative dualizing sheaf
of the family $\left(  \ref{10}\right)  $ and on $\pi_{\ast}\omega
_{\mathcal{E}/\mathfrak{M}_{2}}$ we have a natural metric. The local sections
of $\pi_{\ast}\omega_{\mathcal{E}/\mathfrak{M}_{2}}$ are families of
holomorphic one forms $\omega_{\tau}.$ Then we define $\mathbf{L}^{2}$ metric%
\[
\left\Vert \omega_{\tau}\right\Vert _{\mathbf{L}^{2}}^{2}=\frac{-\sqrt{-1}}{2}%
{\displaystyle\int\limits_{E_{\tau}}}
\omega_{\tau}\wedge\overline{\omega}_{\tau}.
\]
There is a natural compactification $\mathbb{CP}^{1}-\left\{  0,1,\infty
\right\}  $ to $\mathbb{CP}^{1}.$ We can prolong $\pi_{\ast}\omega
_{\mathcal{E}/\mathfrak{M}_{2}}$ to a line bundle $\overline{\pi_{\ast}%
\omega_{\mathcal{E}/\mathfrak{M}_{2}}}$ over $\mathbb{CP}^{1}$ by extendning
the holomorphic sections with finite $\mathbf{L}^{2}$ norm.

It is a well know fact that any elliptic curve $E_{\tau}:=\mathbb{C}%
/\Lambda_{\tau}\,$\ where $\Lambda_{t}$ is the lattice $\left\{
m+n\tau|m,n\in\mathbb{Z},\text{ }\tau\in\mathbb{C}\text{ and }%
\operatorname{Im}\tau>0\right\}  $ can be embedded in $\mathbb{CP}^{2}$ and
the equation in the standard affine open set of $E_{\tau}$ is given by
$y^{2}=4x^{3}-g_{2}(\tau)-g_{3}(\tau),$ where
\[
g_{2}(\tau)=60%
{\displaystyle\sum\limits_{\left(  n,m\right)  \neq(0,0)}}
\frac{1}{\left(  n+m\tau\right)  ^{4}}\text{ and }g_{3}(\tau)=140%
{\displaystyle\sum\limits_{\left(  n,m\right)  \neq(0,0)}}
\frac{1}{\left(  n+m\tau\right)  ^{6}}.
\]
The algebraic discriminant $\Delta(\tau)$ of the elliptic curve $E_{\tau}$ is
defined as the discriminant of the polynomial $4x^{3}-g_{2}(\tau)-g_{3}%
(\tau).$ Thus we get the explicit formula for the algebraic discriminant
\begin{equation}
\Delta(\tau)=g_{2}(\tau)^{3}-27g_{3}(\tau)^{2}. \label{3}%
\end{equation}

The relation between the algebraic discriminants given by $\left(
\ref{2}\right)  $ and $\left(  \ref{3}\right)  $ is given by the following; It
is a well known fact that $\Gamma(2)$ is a normal subgroup in $\mathbb{PSL}%
_{2}\left(  \mathbb{Z}\right)  $ and $\mathbb{PSL}_{2}\left(  \mathbb{Z}%
\right)  \left/  \Gamma(2)\right.  =S_{3},$ where $S_{3}$ is the symmetric
group. Thus we have a finite Galois covering
\begin{equation}
\pi_{2}:\Gamma(2)\left\backslash \mathfrak{h}\right.  \rightarrow
\mathbb{PSL}_{2}\left(  \mathbb{Z}\right)  \left\backslash \mathfrak{h}%
\right.  . \label{4}%
\end{equation}
Thus $\Delta_{alg}(\tau)=\pi_{2}^{\ast}\left(  g_{2}(\tau)^{3}-27g_{3}%
(\tau)^{2}\right)  $ is a section $\Delta$ of $\overline{\pi_{\ast}%
\omega_{\mathcal{E}/\mathfrak{M}_{2}}}$ which vanishes on $\left\{
0,1,\infty\right\}  .$ We will call $\Delta_{alg}(\tau)$ the algebraic
discriminant of the elliptic curve with fixed basis in $H^{1}\left(
E,\mathbb{Z}/2\mathbb{Z}\right)  $.

Let us consider the function $\pi^{\ast}\left(  \Delta\right)  $ on
$\mathfrak{h}.$ Then $\pi^{\ast}\left(  \Delta\right)  $ will be up to a
constant equal to the cusp form of weight $12.$ $\pi^{\ast}\left(
\Delta\right)  $ will be called the analytic discriminant. \ The relations
between the analytic and algebraic discriminants is given by the following
interpretation of the Kronecker limit formula:

\begin{theorem}
\label{kron1}$\det\Delta_{\tau}=\left\Vert \pi^{\ast}\left(  \Delta\right)
\right\Vert _{\mathbf{L}^{2}}^{2}.$
\end{theorem}

\textbf{4. The Analogue of Baily-Borel Compactification of the Moduli Space of
Polarized CY Manifolds}

Baily and Borel constructed a compactification of locally symmetric spaces
quotient by an arithmetic group by using cusp forms, i.e. automorphic forms
which vanish at the cusps. In case of polarized CY manifolds it was proved in
\cite{LSTY} that the completion of the Teichm\"{u}ller space with respect to
the Hodge metric is a domain of holomorphy. We proved in this paper that
sections of the power of the relative dualizing sheaf with finite
$\mathbf{L}^{2}$ norm are the analogue of cusp forms. We show that some power
of the relative dualizing sheaf and its holomorphic sections with finite
$\mathbf{L}^{2}$ norm define a holomorphic embedding of the moduli space
$\mathfrak{M}_{L}\left(  \text{M}\right)  $ of polarized CY manifolds into
$\mathbb{P}^{m_{0}}.$ The projective closure $\overline{\mathfrak{M}%
_{L}\left(  \text{M}\right)  }$ of the image of $\mathfrak{M}_{L}\left(
\text{M}\right)  $ in $\mathbb{P}^{m_{0}}$ will be the analogue of the
Baily-Borel compactification. We also prove the analogue of Borel extension
Theorem, namely any map of $\left(  D^{\ast}\right)  ^{k}\times D^{h^{n-1,1}%
-k}$ into $\mathfrak{M}_{L}\left(  \text{M}\right)  $ can be analytically
prolonged to a holomorphic map of $D^{k}\times D^{h^{n-1,1}-k}$ to the
Baily-Borel compactification $\overline{\mathfrak{M}_{L}\left(  \text{M}%
\right)  }$ of $\mathfrak{M}_{L}\left(  \text{M}\right)  .$ This implies that
the compactification we constructed is minimal model.

\textbf{5. Relations between Algebraic and Analytic Discriminants}

The generalization of the relation between the analytic and algebraic
discriminant on elliptic curves to CY three folds is done in the last section
of this paper and it follows from the proof that the regularized determinants
of the Laplacians of CY metrics acting on $(0,1)$ forms on a CY manifold M are
bounded. The proof of the boundedness is based on the computation of the short
term asymptotic of the trace $Tr(\exp(-t\Delta_{\tau,q}))$ of the heat kernel
of a CY metric for a CY threefold. We established that
\[
Tr(\exp(-t\Delta_{\tau,q}))=\frac{a_{-3}}{t^{3}}+\frac{a_{-2}}{t^{2}}%
+\frac{a_{-1}}{t}+a_{0}+O(t),
\]
where
\[
a_{-3}=%
{\displaystyle\int\limits_{\text{M}}}
L^{3}=vol(g_{\tau}),\text{ }a_{-2}=%
{\displaystyle\int\limits_{\text{M}}}
c_{1}\left(  \text{M}\right)  \wedge L^{2}=0,
\]%
\[
a_{-1}(g)=-\frac{1}{720\pi}\int\limits_{\text{M}}c_{2}(\text{M})\wedge L
\]
and $a_{0}$ are constants. As a consequence of this formula we get that there
exists a non-zero holomorphic section $\eta^{\otimes N}$ of the $N^{th}$ power
$\omega_{\left.  \mathcal{X}\right/  \mathfrak{M}_{L}\left(  \text{M}\right)
}^{\otimes N}$ of the relative dualizing sheaf such that its $\mathbf{L}^{2}$
norm $\left\Vert \eta^{\otimes N}\right\Vert _{\mathbf{L}^{2}}^{2}$ is
$\left(  \det(\Delta_{0,1}(\tau))\right)  ^{N},$ where
\[
N=\#\left(  \left.  \Gamma\right/  \left[  \Gamma,\Gamma\right]  \right)  .
\]
Recall that the moduli space $\mathfrak{M}_{L}\left(  \text{M}\right)  $ of
polarized CY manifolds is obtained from the Teichm\"{u}ller space
$\mathcal{T}\left(  \text{M}\right)  $ by the action of some arithmetic group
$\Gamma$ of rank at least two$.$ From a Theorem proved by Kazhdan it follows
that the abelian group $\left.  \Gamma\right/  \left[  \Gamma,\Gamma\right]  $
is finite.

\begin{conjecture}
\label{to}Suppose that M is a CY manifold of complex dimension $n$ with fixed
polarization class $L$. Suppose that $g$ is a CY metric such that the
cohomology class of the imaginary part is $L.$ Then the coefficients
$a_{-k}(g)$ for $0\leq k\leq n$ of the short term asymptotic expansion of the
trace of the heat kernel are given by the formula:
\[
Tr(k_{t}(x,y))=\frac{a_{-n}(g)}{t^{n}}+...+\frac{a_{-k}(g)}{t^{k}}%
+...+a_{0}(g)+...,
\]
where $a_{-k}(g)=b_{k}\int\limits_{\text{M}}c_{n-k}($M)$\wedge L^{k}$ for
$k=1,...,n,$ $b_{k}$ are some constants which depend on the dimension of M,
and $c_{k}($M) are the Chern classes of M.
\end{conjecture}

Conjecture \ref{to} implies that the regularized determinant $\det
\Delta_{(p,q)}(\tau)$ of the Laplacian of a CY metric with a fixed class of
cohomology of its imaginary part and actining on $(p,q)$ forms is a bounded
function on the moduli space $\mathfrak{M}_{L}\left(  \text{M}\right)  $ of
polarized CY manifolds.

\begin{acknowledgement}
The author express his gratitude to Kefeng Liu for useful remarks,
encouragement and interest to this paper.
\end{acknowledgement}

\section{Basic Definitions and Notions}

\subsection{Kuranishi Space and Flat Local Coordinates}

The following Theorem is proved in \cite{To89}:

\begin{theorem}
\label{tod1}Let M be a CY manifold and let $\left\{  \phi_{i}\right\}  $ for
$i=1,...,N$ be a basis of harmonic $(0,1)$ forms with coefficients in
$T^{1,0},$ i.e. $\left\{  \phi_{i}\right\}  \in\mathbb{H}^{1}($M$,T^{1,0}).$
Then the equation $\overline{\partial}\phi(\tau)=\frac{1}{2}\left[  \phi
(\tau),\phi(\tau)\right]  $ has a solution in the form:
\[
\phi(\tau)=\sum_{i=1}^{N}\phi_{i}\tau^{i}+\sum_{|I_{N}|\geqq2}\phi_{I_{N}}%
\tau^{I_{N}}=
\]%
\begin{equation}
\sum_{i=1}^{N}\phi_{i}\tau^{i}+\frac{1}{2}\overline{\partial}^{\ast}%
G[\phi(\tau^{1},...,\tau^{N}),\phi(\tau^{1},...,\tau^{N})], \label{ANTOD}%
\end{equation}
where $I_{N}=(i_{1},...,i_{N})$\ \ \textit{is a multi-index},
\[
\overline{\partial}^{\ast}\phi(\tau^{1},...,\tau^{N})=0,\text{ }\phi_{I_{N}%
}\lrcorner\omega_{\text{M}}=\partial\psi_{I_{N}}%
\]
\
\[
\phi_{I_{N}}\in C^{\infty}(\text{M},\Omega^{0,1}\otimes T^{1,0}),\tau^{I_{N}%
}=(\tau^{i})^{i_{1}}...(\tau^{N})^{i_{N}}%
\]
\textit{and there exists} $\varepsilon>0$ such that for $\ |\tau
^{i}|<\varepsilon,$ \textit{\ } $\phi(\tau)\in C^{\infty}($M$,\Omega
^{0,1}\otimes T^{1,0}).$ $\left(  See\text{ }\cite{Ti}\text{ }and\text{
}\cite{To89}\right)  .$
\end{theorem}

It is a standard fact from Kodaira-Spencer-Kuranishi deformation theory that
for each $\tau=(\tau^{1},...,\tau^{N})\in\mathcal{K}$ as in Theorem
\ref{tod1}, the Beltrami differential $\phi(\tau^{1},...,\tau^{N})$ defines a
new integrable complex structure on M. This means that the points of
$\mathcal{K},$ define a family of integrable in the sense of
Newlander-Nirenberg operators $\overline{\partial}_{\tau}$ on the $C^{\infty}$
family $\mathcal{K}\times$M$\rightarrow$M. Moreover, it was proved by Kodaira,
Spencer and Kuranishi that over $\mathcal{K}$ there exists a complex analytic
family of CY manifolds $\pi:\mathcal{X\rightarrow K}.$ The family
$\pi:\mathcal{X\rightarrow K}$ is called the Kuranishi family. The operators
$\overline{\partial}_{\tau}$ are defined as follows: Let $\{\mathcal{U}_{i}\}$
be an open covering of M, with a local coordinate system in $\mathcal{U}_{i}$
given by $\{z_{i}^{k}\}$ with $k=1,...,n=$dim$_{\mathbb{C}}$M. Assume that
$\phi(\tau^{1},...,\tau^{N})|_{\mathcal{U}_{i}}$ is given by:
\[
\phi(\tau^{1},...,\tau^{N})=\sum_{j,k=1}^{n}(\phi(\tau^{1},...,\tau
^{N}))_{\overline{j}}^{k}\text{ }d\overline{z}^{j}\otimes\frac{\partial
}{\partial z^{k}}.
\]
Then the operators\textit{\ }$\overline{\partial}_{\tau}$ are given by the
following explicit formulas:\textit{\ \ }\
\begin{equation}
\left(  \overline{\partial}\right)  _{\tau,\overline{j}}=\frac{\overline
{\partial}}{\overline{\partial z^{j}}}-\sum_{k=1}^{n}(\phi(\tau^{1}%
,...,\tau^{N}))_{\overline{j}}^{k}\frac{\partial}{\partial z^{k}}.
\label{eq:dibar}%
\end{equation}

\begin{definition}
\label{flat}The coordinates $\tau=(\tau^{1},...,\tau^{N})$ defined in Theorem
\ref{tod1} will be fixed from now on and will be called the flat coordinate
system in $\mathcal{K}$.
\end{definition}

\subsection{The Definition of Hodge and Weil-Petersson Metrics}

It is a well-known fact from Kodaira-Spencer-Kuranishi theory that the tangent
space $T_{\tau,\mathcal{K}\text{ }}$at a point $\tau\in\mathcal{K}$ can be
identified with the space of harmonic (0,1) forms with values in the
holomorphic vector fields, which we will denote by $\mathbb{H}^{1}($M$_{\tau
},T$). We will view each element $\phi\in\mathbb{H}^{1}($M$_{\tau},T$) as a
point wise linear map from $\phi:\Omega_{\text{M}_{\tau}}^{(1,0)}%
\rightarrow\Omega_{\text{M}_{\tau}}^{(0,1)}.$ Given $\phi_{1}$ and $\phi
_{2}\in\mathbb{H}^{1}($M$_{\tau},T$)$,$ the trace of the map $\phi_{1}%
\circ\overline{\phi_{2}}:\Omega_{\text{M}_{\tau}}^{(0,1)}\rightarrow
\Omega_{\text{M}_{\tau}}^{(0,1)}$ at the point $m\in$M$_{\tau}$ with respect
to the metric g is simply:\textit{\ \ }\
\begin{equation}
Tr(\phi_{1}\circ\overline{\phi_{2}})(m)=\sum_{k,l,m=1}^{n}(\phi_{1}%
)_{\overline{l}}^{k}(\overline{\phi_{2})_{\overline{k}}^{m}}g^{\overline{l}%
,k}g_{k,\overline{m}} \label{wp}%
\end{equation}
We will define the Weil-Petersson metric on $\mathcal{K}$ as:
\begin{equation}
\left\langle \phi_{1},\phi_{2}\right\rangle =\int\limits_{\text{M}}Tr(\phi
_{1}\circ\overline{\phi_{2}})vol(g). \label{wp1}%
\end{equation}

\begin{definition}
\label{km}Let $\mathbb{G}$ be a semi-simple Lie group, $\mathbb{K}$ be a
maximal compact group and $\mathbb{K}_{1}$ be a proper subgroup in
$\mathbb{K}$. Let us consider the homogeneous space $\mathbb{G}\left/
\mathbb{K}_{1}\right.  $ and its projection: $\mathbb{P}\mathbf{r}%
:\mathbb{G}\left/  \mathbb{K}_{1}\right.  \rightarrow\mathbb{G}\left/
\mathbb{K}\right.  .$ Let us consider the Cartan decomposition of the Lie
algebra $\mathcal{G}$ of $\mathbb{G}:$%
\begin{equation}
\mathcal{G}=\mathcal{K}\oplus\mathcal{P}, \label{cd}%
\end{equation}
where $\mathcal{K}$ is the Lie algebra of $\mathbb{K}$ and $\left(
\ref{cd}\right)  $ is orthogonal decomposition of the Lie algebra
$\mathcal{G}$ with respect to the Killing form $K(x,y)$. Then the Killing form
is non degenerate negative bilinear on $\mathcal{K}$ and positive
non-degenerate form on $\mathcal{P}.$ The tangent space $T_{id,\mathbb{G}%
\left/  \mathbb{K}_{1}\right.  }$ is isomorphic to $\mathcal{K}\left/
\mathcal{K}_{1}\right.  \oplus\mathcal{P},$ where $\mathcal{K}_{1}$ is the Lie
algebra of $\mathbb{K}_{1}.$ Then if $\left(  \alpha,\beta\right)
\in\mathcal{K}_{1}^{\bot}\oplus\mathcal{P},$ where $\mathcal{K}_{1}^{\bot
}\subset\mathcal{K}$ is the perpenicular space to $\mathcal{K}_{1}$ in
$\mathcal{K},$ we define the Killing norm of $(\alpha,\beta)$ as follows:%
\begin{equation}
\left\Vert \left(  \alpha,\beta\right)  \right\Vert _{K}^{2}=-K(\alpha
,a)+K(\beta,\beta). \label{Hm}%
\end{equation}
Thus $\left(  \ref{Hm}\right)  $ define an invariant metric on $\mathbb{G}%
\left/  \mathbb{K}_{1}\right.  .$ We will call this metric the Hodge metric on
$\mathbb{G}\left/  \mathbb{K}_{1}\right.  .$
\end{definition}

\begin{remark}
\label{HMr}It is easy to see that the Hodge metric on $\mathbb{G}\left/
\mathbb{K}_{1}\right.  $ is a complete metric.
\end{remark}

\begin{definition}
\label{HM}It is a well known and easy fact that the moduli space of Variations
of Hodge Structures of given weight is isomorphic $\mathbb{G}\left/
\mathbb{K}_{1}\right.  ,$ where $\mathbb{G}$ is a semi-simple Lie group and
$\mathbb{K}_{1}$ is a compact subgroup. The period map:%
\[
p:\mathfrak{M}_{L}\left(  \text{M}\right)  \rightarrow\mathbb{G}\left/
\mathbb{K}_{1}\right.
\]
is a well defined and it is a local isomorphism when M is a CY manifold. Then
the restriction of the Hodge metric on $p\left(  \mathfrak{M}_{L}\left(
\text{M}\right)  \right)  $ defines the Hodge metric on the moduli space
$\mathfrak{M}_{L}\left(  \text{M}\right)  $ of polarized CY manifolds.
\end{definition}

In \cite{LSTY} the following Theorem was proved:

\begin{theorem}
\label{NC}The Hodge metric on the moduli space $\mathfrak{M}_{L}\left(
\text{M}\right)  $ of polarized CY manifolds has non-positive curvature and
bounded from above holomorphic sectional curvature by a negative constant.
\end{theorem}

\subsection{Review of the Results in \cite{LTYZ}}

\begin{definition}
\label{Teich}We will define the Teichm\"{u}ller space $\mathcal{T}$(M) of a CY
manifold M as follows: $\mathcal{T}($M$):=\mathcal{I}($M$)/Diff_{0}($M$),$
\textit{where}\
\[
\mathcal{I}(\text{M}):=\left\{  \text{all integrable complex structures on
M}\right\}
\]
\textit{and } Diff$_{0}$(M) \textit{is the group of diffeomorphisms isotopic
to identity. The action of the group Diff(M}$_{0})$ \textit{is defined as
follows; Let }$\phi\in$Diff$_{0}$(M) \textit{then }$\phi$ \textit{acts on
integrable complex structures on M by pull back, i.e. if }
\[
I\in C^{\infty}(\text{M},Hom(T(\text{M}),T(\text{M})),
\]
\textit{then we define } $\phi(I_{\tau})=\phi^{\ast}(I_{\tau}).$
\end{definition}

We will call a pair (M; $\gamma_{1},...,\gamma_{b_{n}}$) a marked CY manifold
where M is a CY manifold and $\{\gamma_{1},...,\gamma_{b_{n}}\}$ is a basis of
$H_{n}$(M,$\mathbb{Z}$)/Tor.

\begin{remark}
\label{mark}Let $\mathcal{K}$ be the Kuranishi space. It is easy to see that
if we choose a basis of $H_{n}$(M,$\mathbb{Z}$)/Tor in one of the fibres of
the Kuranishi family $\mathcal{M\rightarrow K}$ then all the fibres will be
marked, since as a $C^{\infty}$ manifold $\mathcal{X}_{\mathcal{K}}\approxeq
$M$\times\mathcal{K}$.
\end{remark}

\begin{theorem}
\label{teich}There exists a family of marked polarized CY manifolds
\begin{equation}
\mathcal{Z}_{L}\mathcal{\rightarrow}\widetilde{T}(\text{M}), \label{fam2}%
\end{equation}
which possesses the following properties: \textbf{a)} It is effectively
parametrized, \textbf{b) }For any marked CY manifold M of fixed topological
type for which the polarization class $L$ defines an imbedding into a
projective space $\mathbb{CP}^{N},$ there exists an isomorphism of it (as a
marked CY manifold) with a fibre M$_{s}$ of the family $\mathcal{Z}_{L}.$
\textbf{c) }The base has dimension $h^{n-1,1}.$
\end{theorem}

\begin{corollary}
\label{teich1}Let $\mathcal{Y\rightarrow}$X be any family of marked CY
manifolds, then there exists a unique holomorphic map $\phi:$X$\rightarrow
\widetilde{T}($M$)$ up to a biholomorphic map $\psi$ of M which induces the
identity map on $H_{n}($M$,\mathbb{Z}).$
\end{corollary}

From now on we will denote by $\mathcal{T}$(M) the irreducible component of
the Teichm\"{u}ller space that contains our fixed CY manifold M.

\subsection{Construction of the Moduli Space of Polarized CY Manifolds}

\begin{definition}
We will define the mapping class group $\Gamma_{1}$ of any compact C$^{\infty
}$ manifold M as follows: $\Gamma_{1}=Diff_{+}($M$)/Diff_{0}($M$),$ where
$Diff_{+}($M$)$ is the group of diffeomorphisms of M preserving the
orientation of M and $Diff_{0}($M$)$ is the group of diffeomorphisms isotopic
to identity.
\end{definition}

\begin{definition}
Let $L\in H^{2}($M$,\mathbb{Z})$ be the imaginary part of a K\"{a}hler metric.
Let $\Gamma_{L}:=\{\phi\in\Gamma_{1}|\phi(L)=L\}.$
\end{definition}

It is a well know fact that the moduli space of polarized algebraic manifolds
$\mathcal{M}_{L}($M$)=\Gamma_{L}\left\backslash \mathcal{T}(M)\right.  .$

\begin{theorem}
\label{Vie}There exists a subgroup of finite index $\Gamma$ of $\ \Gamma_{L}$
such that $\Gamma$ acts freely on $\mathcal{T}$(M) and $\Gamma\backslash
\mathcal{T}\left(  \text{M}\right)  =\mathfrak{M}_{L}\left(  \text{M}\right)
$ is a non-singular quasi-projective variety.
\end{theorem}

\begin{remark}
\label{Vie1}Theorem \ref{Vie} implies that we constructed a family of
non-singular CY manifolds $\pi:\mathcal{X\rightarrow}\mathfrak{M}_{L}\left(
\text{M}\right)  $ over a quasi-projective non-singular variety $\mathcal{M}%
$(M). Moreover it is easy to see that $\mathcal{X\subset}\mathbb{CP}^{N}%
\times\mathfrak{M}_{L}\left(  \text{M}\right)  .$ \textit{So} $\mathcal{X}$
\ \textit{is also quasi-projective. From now on we will work only with this
family.}
\end{remark}

\subsection{Metrics on Vector Bundles with Logarithmic Growth}

In Theorem \ref{Vie} we constructed the moduli space $\mathfrak{M}_{L}\left(
\text{M}\right)  $ of CY manifolds. From the results in \cite{W} and Theorem
\ref{Vie} we know that $\mathfrak{M}_{L}\left(  \text{M}\right)  $ is a
quasi-projective non-singular variety. Using Hironaka's resolution theorem, we
may suppose that $\mathfrak{M}_{L}\left(  \text{M}\right)  \subset
\overline{\mathfrak{M}_{L}\left(  \text{M}\right)  },$ where $\overline
{\mathfrak{M}_{L}\left(  \text{M}\right)  }-\mathfrak{M}_{L}\left(
\text{M}\right)  =\mathfrak{D}$ is a divisor with normal crossings. We need
now to show how we will extend the determinant line bundle $\mathcal{L}$ to a
line bundle $\overline{\mathcal{L}\text{ }}$ to $\overline{\mathfrak{M}%
_{L}\left(  \text{M}\right)  }.$ For this reason we are going to recall the
following definitions and results from \cite{Mum}. We will look at polydisks
D$^{N}\subset\overline{\mathfrak{M}_{L}\left(  \text{M}\right)  },$ where D is
the unit disk, $N=\dim\overline{\mathfrak{M}_{L}\left(  \text{M}\right)  }$
and such that
\[
D^{N}\cap\mathfrak{D}_{\infty}=\{\text{union of hyperplanes};\tau
_{1}=0,...,\tau_{k}=0\}.
\]
Hence, $D^{N}\cap\mathcal{M}(M)=(D^{\ast})^{k}\times D^{N-k}.$ On $D^{\ast}$
we have the Poincare metric
\[
ds^{2}=\frac{\left\vert dz\right\vert ^{2}}{\left\vert z\right\vert
^{2}\left(  \log\left\vert z\right\vert \right)  ^{2}}%
\]
and on D we have the simple metric $\left\vert dz\right\vert ^{2},$ giving us
a product metric on ($D^{\ast}$)$^{k}\times D^{N-k}$ which we call
$\omega^{(P)}.$

A complex-valued C$^{\infty}$ p-form $\eta$ on $\mathcal{M}$(M) is said to
have Poincare growth on $\overline{\mathfrak{M}_{L}\left(  \text{M}\right)
}-\mathfrak{M}_{L}\left(  \text{M}\right)  $ if there is a set of \textit{if
}polydisks\textit{ \ } $\mathcal{U}_{\alpha}\subset\overline{\mathfrak{M}%
_{L}\left(  \text{M}\right)  }$ covering $\overline{\mathfrak{M}_{L}\left(
\text{M}\right)  }-\mathfrak{M}_{L}\left(  \text{M}\right)  $ such that in
each \textit{\ }$\mathcal{U}_{\alpha}$ an estimate of the following type
holds:
\[
\left\vert \eta(\tau_{1},...,\tau_{N}\right\vert \leq C_{\alpha}%
\omega_{\mathcal{U}_{\alpha}}^{(P)}(\tau_{1},\overline{\tau_{1}})..._{\alpha
}\omega_{\mathcal{U}_{\alpha}}^{(P)}(\tau_{N},\overline{\tau_{N}}).
\]
This property is independent of the covering $\mathcal{U}_{\alpha}$ of
$\overline{\mathfrak{M}_{L}\left(  \text{M}\right)  }-\mathfrak{M}_{L}\left(
\text{M}\right)  $ but depends on the compactification $\overline
{\mathfrak{M}_{L}\left(  \text{M}\right)  }.$ If $\eta_{1}$ and $\eta_{2}$
both have Poincare growth on $\overline{\mathfrak{M}_{L}\left(  \text{M}%
\right)  }-\mathfrak{M}_{L}\left(  \text{M}\right)  ,$ then so does $\eta
_{1}\wedge\eta_{2}$. The basic property of the Poincare growth is the following:

\begin{theorem}
\label{Mum}A p-form $\eta$ with a Poincare growth on $\overline{\mathfrak{M}%
_{L}\left(  \text{M}\right)  }-\mathfrak{M}_{L}\left(  \text{M}\right)
=\mathfrak{D}$ \textit{has the property that for every C}$^{\infty}$ (r-p)
\textit{form} $\psi$ \textit{on \ }$\overline{\mathfrak{M}_{L}\left(
\text{M}\right)  }$ \textit{we have:}
\[
\int_{\overline{\mathfrak{M}_{L}\left(  \text{M}\right)  }-\mathfrak{M}%
_{L}\left(  \text{M}\right)  }\left\vert \eta\wedge\psi\right\vert <\infty.
\]
\textit{Hence, }$\eta$ \textit{defines a current [}$\eta$] \textit{on
}$\overline{\mathfrak{M}_{L}\left(  \text{M}\right)  }.$
\end{theorem}

\textbf{Proof:}For the proof see \cite{Mum}. $\blacksquare.$

A complex valued C$^{\infty}$ p-form\textit{\ }$\eta$ on\textit{\ }%
$\overline{\mathfrak{M}_{L}\left(  \text{M}\right)  }$\textit{\ is good on M
if both }$\eta$\textit{\ }and $d\eta$ have Poincare growth\textit{. }Let
$\mathcal{E}$ be a vector bundle on $\mathcal{M}$(M) with a Hermitian metric
h. We will call h a good metric on $\overline{\mathfrak{M}_{L}\left(
\text{M}\right)  }$ if the following holds:

\begin{enumerate}
\item \ If\textbf{ }for all x$\in\overline{\mathfrak{M}_{L}\left(
\text{M}\right)  }-\mathfrak{M}_{L}\left(  \text{M}\right)  ,$ there exists
sections $e_{1},...,e_{m}$ of $\mathcal{E}$ \ which form a basis of
\textit{\ \ }$\mathcal{E}\left\vert _{D^{r}-(D^{r}\cap\mathfrak{D}_{\infty}%
)}\right.  .$

\item In a neighborhood \textit{D}$^{r}$ of x in which $\overline
{\mathfrak{M}_{L}\left(  \text{M}\right)  }-\mathfrak{M}_{L}\left(
\text{M}\right)  $ is given by\textit{\ }
\[
z_{1}\times...\times z_{k}=0.
\]

\item The metric \textit{\ }h$_{i\overline{j}}=$h($e_{i},e_{j}$) has the
following properties: \textbf{a.}
\[
\left\vert h_{i\overline{j}}\right\vert \leq C\left(  \sum_{i=1}^{k}%
\log\left\vert z_{i}\right\vert \right)  ^{2m},\text{ }\left(  \det\left(
h\right)  \right)  ^{-1}\leq C\left(  \sum_{i=1}^{k}\log\left\vert
z_{i}\right\vert \right)  ^{2m}%
\]
for some $C>0,$ $m\geq0.$ \textbf{b.}\ The 1-forms\textit{\ }$\left(  \left(
dh\right)  h^{-1}\right)  _{i\overline{j}}$ are good forms on \textit{\ }%
$\overline{\mathfrak{M}_{L}\left(  \text{M}\right)  }\cap D^{N}.$
\end{enumerate}

It is easy to prove that there exists a unique extension $\overline
{\mathcal{E}}$ of $\mathcal{E}$ \ on $\overline{\mathfrak{M}_{L}\left(
\text{M}\right)  },$ i.e. $\overline{\mathcal{E}}$ is defined locally
as\ holomorphic sections of $\mathcal{E}$ which have a finite norm in h.

\begin{theorem}
\label{Mum100}Let ($\mathcal{E}$,h) be a vector bundle with a good metric on
$\mathfrak{M}_{L}\left(  \text{M}\right)  $, then the Chern classes c$_{k}%
$($\mathcal{E}$,h) are good forms on $\overline{\mathfrak{M}_{L}\left(
\text{M}\right)  }$ \textit{and the currents }$\left[  c_{k}(\mathcal{E}%
,\mathbf{L}^{2})\right]  $ \textit{represent the cohomology classes }
\[
c_{k}(\mathcal{E},\mathbf{L}^{2})\in H^{2k}(\overline{\mathfrak{M}_{L}\left(
\text{M}\right)  },\mathbb{Z)}.
\]

\end{theorem}

\textbf{Proof: }For the proof see \cite{Mum}. $\blacksquare.$

\subsection{Applications of Mumford's Results to the Moduli of CY}

In \cite{ls2} and \cite{To04} the following Theorem was proved:

\begin{theorem}
\label{Nik}Let $\pi:\mathcal{X\rightarrow}\mathfrak{M}_{L}\left(
\text{M}\right)  $ be the flat family of non-singular CY manifolds. Let the
relative dualizing sheaf $\omega_{\mathcal{X}\text{/}\mathcal{M}\text{(M)}%
}:=\pi_{\ast}\Omega_{\mathcal{X}\text{/}\mathcal{M}\text{(M)}}^{n,0}$
\textit{be equipped with the metric }$\mathbf{L}^{2}$\textit{ defined by}%
\begin{equation}
\left\Vert \omega\right\Vert _{\mathbf{L}^{2}}^{2}:=(-1)^{\frac{n(n-1)}{2}%
}\left(  \frac{\sqrt{-1}}{2}\right)  ^{n}%
{\displaystyle\int\limits_{\text{M}}}
\omega\wedge\overline{\omega} \label{l2}%
\end{equation}
\textit{ Then }$\mathbf{L}^{2}$\textit{\ is a good metric. }
\end{theorem}

\section{The Analogue of Baily-Borel Compactification of the Moduli Space of
CY Manifolds}

\subsection{Construction of the Analogue of the Dedekind $\eta$ Function for
CY Manifolds}

\begin{theorem}
\label{Nik2}Let M be a CY manifold. Let $N=\#\Gamma/[\Gamma,\Gamma].$ Then
$\omega_{\mathcal{X}/\mathfrak{M}_{L}(\text{M})}^{\otimes N}$ is a trivial
complex analytic line bundle over $\mathfrak{M}_{L}\left(  \text{M}\right)  $.
\end{theorem}

\textbf{Proof: }The proof of Theorem \ref{Nik2} is based on the following
Theorem proved in \cite{To89}:

\begin{theorem}
\label{to89}The Chern class of the relative dualizing sheaf $\omega
_{\mathcal{X}/\mathfrak{M}_{L}(\text{M})}$ is the imaginary part of the
Weil-Petersson metric on $\mathfrak{M}_{L}($M$).$
\end{theorem}

According to the results proved in \cite{BT} we have
\begin{equation}
dd^{c}\log\left(  \det\Delta_{(0,1)}(\tau)\right)  =\operatorname{Im}%
W-P.\label{a}%
\end{equation}
So from $\left(  \ref{a}\right)  $, the fact that the $\mathbf{L}^{2}$ metric
is good, Theorem \ref{Mum100} and Theorem \ref{to89} we get that the Chern
class of the relative dualizing sheaf $\omega_{\mathcal{X}/\mathfrak{M}%
_{L}(\text{M})}$ is zero in $H^{2}\left(  \mathfrak{M}_{L}(\text{M}%
),\mathbb{Z}\right)  .$ This means that $\omega_{\mathcal{X}/\mathfrak{M}%
_{L}(\text{M})}$ is a trivial $C^{\infty}$ line bundle on $\mathfrak{M}_{L}%
($M$).$ Let us denote by $\sigma:\mathcal{T}($M$)\rightarrow\mathfrak{M}%
_{L}\left(  \text{M}\right)  =\Gamma\left\backslash \mathcal{T}(\text{M}%
)\right.  $ the natural projection map. So the line bundle $\sigma^{\ast
}\left(  \omega_{\mathcal{X}/\mathfrak{M}_{L}(\text{M})}\right)  $ will be
trivial on $\mathcal{T}($M$)$, i.e. $\sigma^{\ast}\left(  \omega
_{\mathcal{X}/\mathfrak{M}_{L}(\text{M})}\right)  \mathbb{\approxeq
}\mathcal{T}($M$)\times\mathbb{C}$ and
\begin{equation}
\omega_{\mathcal{X}/\mathfrak{M}_{L}(\text{M})}\approxeq\Gamma\left\backslash
\mathbb{C\times}\mathcal{T}(\text{M})\right.  ,\label{n}%
\end{equation}
where $\Gamma$ acts in a natural way on the Teichm\"{u}ller space and it acts
by a character
\[
\chi\in Hom(\Gamma,\mathbb{C}_{1}^{\ast})\approxeq Hom(\Gamma/[\Gamma
,\Gamma],\mathbb{C}_{1}^{\ast})
\]
of the group $\Gamma$ on the fibre $\mathbb{C}.$ A Theorem of Kazhdan states
that $\Gamma/[\Gamma,\Gamma]$ is a finite group if the rank of $\Gamma$ is
bigger or equal to 2. For CY manifolds $\Gamma$ is an arithmetic group of rank
$\geq2$ according to \cite{Sul}. From here we deduce that $\omega
_{\mathcal{X}/\mathfrak{M}_{L}(\text{M})}^{\otimes N}$ will be a trivial
complex analytic bundle on $\mathfrak{M}_{L}($M$)$, where $N=\#\left.
\Gamma\right/  [\Gamma,\Gamma].$ Theorem \ref{Nik2} is proved. $\blacksquare$

\begin{theorem}
\label{sec}Let M be a CY manifold. Let $\mathfrak{M}_{L}($M) be the moduli
space of polarized CY manifolds such that $\tau_{0}\in$ $\mathfrak{M}_{L}($M)
corresponds to M. Let $\omega_{\mathcal{X}/\mathfrak{M}_{L}(\text{M})}$ be the
relative dualizing sheaf of the family $\mathcal{X}\rightarrow\mathfrak{M}%
_{L}($M$).$ Then there exists a non zero section $\eta^{N}\in H^{0}\left(
\mathfrak{M}_{L}(\text{M})\text{,}\omega_{\mathcal{X}/\mathfrak{M}%
_{L}(\text{M})}^{\otimes N}\right)  $ such that $\left\Vert \eta
(\tau)\right\Vert _{\mathbf{L}^{2}}^{2}=\det\Delta_{(0,1)}(\tau).$
\end{theorem}

\textbf{Proof: }Let $\{\phi_{i}\}$ be a basis of harmonic Dolbault
representatives with respect to the CY metric corresponding to the
polarization class $L$ of $\mathbb{H}^{1}\left(  \text{M,}T_{\text{M}}%
^{1,0}\right)  $. As it was proved in \cite{To89} $\{\phi_{i}\}$ defines a
coordinate system $(\tau^{1},...,\tau^{N})$ in the local deformation space
$\mathcal{K}$ which we will call the Kuranishi space. We know that the CY
metric with a fixed polarization class depends real analytically on the
coordinates $(\tau^{1},...,\tau^{N},\overline{\tau^{1}},...,\overline{\tau
^{N}}).$ From here it follows that the regularized determinants $\det
\Delta_{(0,q)}(\tau)$ depend also real analytically on the coordinates
$(\tau^{1},...,\tau^{N},\overline{\tau^{1}},...,\overline{\tau^{N}}).$ The
main result in \cite{BT} is the following Theorem:

\begin{theorem}
\label{in}We have: $\left.  \frac{\partial^{2}}{\overline{\partial\tau^{j}%
}\partial\tau^{i}}\log\left(  \det\Delta_{\tau,\left(  0,1\right)  }\right)
\right\vert _{\tau=0}=\left\langle \phi_{i},\phi_{j}\right\rangle .$
\end{theorem}

In \cite{To89} we proved the following Theorem:

\begin{theorem}
\label{in1}We have$\left.  \frac{\partial^{2}}{\overline{\partial\tau^{j}%
}\partial\tau^{i}}\left(  \log\left\langle \omega_{\tau},\omega_{\tau
}\right\rangle \right)  \right\vert _{\tau=0}=\left\langle \phi_{i},\phi
_{j}\right\rangle .$
\end{theorem}

Theorems \ref{in} and \ref{in1} imply\ that for each point $\tau
\in\mathfrak{M}_{L}($M$)$ there exists an open set $\mathcal{U}_{\tau}$ such
that for $\tau\in\mathcal{U}_{\tau}$ and
\begin{equation}
\left\Vert \eta(\tau)\right\Vert _{\mathbf{L}^{2}}^{2}=\frac{\left\vert
f_{\mathcal{U}_{\tau}}(\tau)\right\vert ^{2}}{\left\langle \omega_{\tau
},\omega_{\tau}\right\rangle }=\frac{\left\vert f_{\mathcal{U}_{\tau}}%
(\tau)\right\vert ^{2}}{\left\Vert \omega_{\tau}\right\Vert _{\mathbf{L}^{2}%
}^{2}},\label{NIK4}%
\end{equation}
where $f_{\mathcal{U}_{\tau}}(\tau)$ is a holomorphic function in
$\mathcal{U}_{\tau}$. Let $\left\{  U_{\alpha}\right\}  $ be a covering of
$\mathfrak{M}_{L}\left(  \text{M}\right)  $ by polydisks. Then $\left(
\ref{NIK4}\right)  $ implies that
\begin{equation}
\left.  \left\Vert \eta(\tau)\right\Vert _{\mathbf{L}^{2}}^{2}\right\vert
_{U_{\alpha}}=\frac{\left\vert f_{\mathcal{\alpha}}(\tau)\right\vert ^{2}%
}{\left\langle \omega_{\tau},\omega_{\tau}\right\rangle }=\frac{\left\vert
f_{\mathcal{\alpha}}(\tau)\right\vert ^{2}}{\left\Vert \omega_{\tau
}\right\Vert _{\mathbf{L}^{2}}^{2}}.\label{NIK4a}%
\end{equation}
Theorem \ref{Nik2} and $\left(  \ref{NIK4a}\right)  $ imply that there exists
a global section $\eta^{N}\in H^{0}\left(  \text{M,}\omega_{\mathcal{X}%
/\mathfrak{M}_{L}(\text{M})}^{\otimes N}\right)  $ such that $\left.  \eta
^{N}\right\vert _{U_{\alpha}}=f_{\alpha}^{N}(\tau).$ Then $\left(
\ref{NIK4a}\right)  $ implies Theorem \ref{sec}. $\blacksquare$

We proved the following Theorem in \cite{LSTY}:

\begin{theorem}
\label{dht}Let $\widetilde{\mathcal{T}(\text{M})}$ be the completion of the
Teichm\"{u}ller space $\mathcal{T}($M$)$ with respect to the Hodge metric.
Then $\widetilde{\mathcal{T}(\text{M})}$ is a domain of holomorphy,
$\mathcal{T}($M$)$ is an open and everywhere dense subset in $\widetilde
{\mathcal{T}(\text{M})}$ and $\widetilde{\mathcal{T}(\text{M})}-\mathcal{T}%
($M$)$ is a countable union of complex subspaces$.$
\end{theorem}

\begin{definition}
\label{Disc}The arithmetic group $\Gamma$ acts on the completion
$\widetilde{\mathcal{T}(\text{M})}$ of the Teichm\"{u}ller space
$\mathcal{T}($M$).$ Let us denote by $\widetilde{\mathfrak{M}_{L}(\text{M}%
)}:=\Gamma\left\backslash \widetilde{\mathcal{T}(\text{M})}\right.  .$ Let
$\overline{\mathfrak{M}_{L}(\text{M})}$ be the projective compactification of
$\mathfrak{M}_{L}($M$)$ such that
\[
\overline{\mathfrak{M}_{L}(\text{M})}-\mathfrak{M}_{L}(\text{M})=\mathfrak{D}%
\]
is a divisor of normal crossings. We will define $\mathfrak{D}_{\infty}$ as
follows; A point $\tau_{\infty}\in\mathfrak{D}_{\infty}$ if and only if
\textbf{1. }around $\tau_{\infty}$ we can find a disk $\mathcal{D}$ such that
$\ $%
\[
\tau_{\infty}\in\mathcal{D},\text{ }\mathcal{D}-\tau_{\infty}\subset
\mathfrak{M}_{L}\left(  \text{M}\right)
\]
over $\mathcal{D}$ $-\tau_{\infty}$ the restriction of the family
$\mathcal{X}\rightarrow\mathfrak{M}_{L}\left(  \text{M}\right)  $ of the
family of polarized CY manifolds on $\mathcal{D}-\tau_{\infty}$ has a
monodromy group of infinite order in $H_{n}($M$_{\tau},\mathbb{Q})$ or
\textbf{2. }$\mathfrak{D}_{f}$ is the codimension one component of
$\widetilde{\mathfrak{M}_{L}(\text{M})}-\mathfrak{M}_{L}($M$)$ that is fixed
point set by a subroup $G\subset\Gamma.$
\end{definition}

\begin{theorem}
\label{Nik3}Let M be a CY manifold$.$ Then there exists a holomorphic section
$\eta^{N}$ of $\omega_{\mathcal{X}/\mathfrak{M}_{L}(\text{M})}^{\otimes
N\text{ }}$ on $\mathfrak{M}_{L}($M$)$ such that it can be prolonged to a
holomorphic section $\overline{\eta^{N}}$ of the line bundle $\overline
{\omega_{\mathcal{X}/\mathfrak{M}_{L}(\text{M})}^{\otimes N\text{ }}}$ such
that for each point $m\in\mathfrak{M}_{L}\left(  \text{M}\right)  ,$ $\eta
^{N}(m)\neq0,$ and the support of the zero set of $\overline{\eta^{N}}$ is the
divisor $\mathfrak{D}_{\infty}.$
\end{theorem}

\textbf{Proof:} The proof of Theorem \ref{Nik3} is based on the following Lemmas:

\begin{lemma}
\label{Nik0}$\overline{\omega_{\mathcal{X}/\mathfrak{M}_{L}(\text{M}%
)}^{\otimes N\text{ }}}\approxeq\mathcal{O}_{\overline{\mathfrak{M}%
_{L}(\text{M})}}\left(
{\displaystyle\sum\limits_{j}}
k_{j}D_{j}\right)  ,$ where $k_{j}\geq0$ and $D_{j}$ are components of
$\overline{\mathfrak{M}_{L}(\text{M})}-\mathfrak{M}_{L}($M$).$
\end{lemma}

\textbf{Proof: }Let $\mathfrak{D}=\bigcup\limits_{i}D_{i}$ be the
decomposition of the divisor $\mathfrak{D}$ on irreducible components on
$\overline{\mathfrak{M}_{L}(\text{M})}$. Theorem \ref{Nik2} implies that the
line bundle $\omega_{\mathcal{X}/\mathfrak{M}_{L}(\text{M})}^{\otimes N}$ is
holomorphic trivial bundle on $\mathfrak{M}_{L}\left(  \text{M}\right)  $ and
$N=\#\left(  \Gamma/\left[  \Gamma,\Gamma\right]  \right)  .$ So we can
conclude that
\begin{equation}
\overline{\omega_{\mathcal{X}/\mathfrak{M}_{L}(\text{M})}^{\otimes N}%
}\approxeq\mathcal{O}_{\overline{\mathfrak{M}_{L}(\text{M})}}\left(  \sum
_{j}k_{j}D_{j}\right)  , \label{NIK2}%
\end{equation}
where $D_{j}$ are the components of $\mathfrak{D}$. We will prove that the
multiplicities $k_{i}$ are non negative integers. Indeed we know from Theorem
\ref{Nik} that the $\mathbf{L}^{2}$ metric defined on the line bundle
$\omega_{\mathcal{X}/\mathfrak{M}_{L}(\text{M})}$ is a good one in the sense
of Mumford. So the Chern form $c_{1}(\omega_{\mathcal{X}/\mathfrak{M}%
_{L}(\text{M})},\mathbf{L}^{2})$ of the good metric h defined by $\left(
\ref{l2}\right)  $ is a positive current on $\overline{\mathfrak{M}_{L}\left(
\text{M}\right)  }$. The Poincare dual of the cohomology of the current
\[
\left[  c_{1}\left(  \omega_{\mathcal{X}/\mathfrak{M}_{L}(\text{M}%
)},\mathbf{L}^{2}\right)  \right]  \in H^{2}\left(  \overline{\mathfrak{M}%
_{L}\left(  \text{M}\right)  },\mathbb{Z}\right)
\]
is
\begin{equation}
\mathcal{P}\left(  \left[  c_{1}(\left(  \omega_{\mathcal{X}/\mathfrak{M}%
_{L}(\text{M})},\mathbf{L}^{2}\right)  \right]  \right)  =\sum_{j}k_{j}%
[D_{j}]\in H_{2n-2}\left(  \overline{\mathfrak{M}_{L}\left(  \text{M}\right)
},\mathbb{Z}\right)  . \label{NIK3}%
\end{equation}
where the coefficients $k_{i}$ are defined as in $\left(  \ref{NIK2}\right)
$. The positivity of the current $c_{1}\left(  \omega_{\mathcal{X}%
/\mathfrak{M}_{L}(\text{M})},\mathbf{L}^{2}\right)  $ implies that its
Poincare dual current $\sum_{j}k_{j}[D_{j}]$ is positive. From here we can
conclude that the coefficients $k_{i}$ are positive integers. Indeed, let
$[\omega_{D_{i}}]\in H^{2n-2}($M,$\mathbb{Z})$ be such classes of cohomology
that:
\begin{equation}
\int\limits_{D_{j}}[\omega_{D_{i}}]=\delta_{ij}. \label{NIK6}%
\end{equation}
Since the current $\sum_{j}k_{i}[D_{j}]$ is positive $\left(  \ref{NIK6}%
\right)  $ implies
\begin{equation}
\left\langle \sum_{j}k_{i}[D_{j}],[\omega_{D_{i}}]\right\rangle =k_{i}\geq0.
\label{NIK6a}%
\end{equation}
Lemma \ref{Nik0} is proved. $\blacksquare$

\begin{corollary}
\label{nik}There exists a section $\eta^{N\text{ }}$of $\overline
{\omega_{\mathcal{X}/\mathfrak{M}_{L}(\text{M})}^{\otimes N}}$ over
$\overline{\mathfrak{M}_{L}(\text{M})}$ such that it vanishes on components of
$\overline{\mathfrak{M}_{L}(\text{M})}-\mathfrak{M}_{L}($M$).$
\end{corollary}

\begin{lemma}
\label{NIK-}The zero set of the section $\overline{\eta^{N}}$ constructed in
Corollary \ref{nik}\ of the line bundle $\overline{\omega_{\mathcal{X}%
/\mathfrak{M}_{L}(\text{M})}^{\otimes N}}^{\text{ }}$ is a non zero divisor
with the same support as $\mathfrak{D}_{\infty}.$
\end{lemma}

\textbf{Proof: }The proof of Lemma \ref{NIK-} is based on the followig two Propositions:

\begin{proposition}
\label{NikA}Let M be a CY manifold. Let $\mathcal{U}_{\tau_{\infty}}$ be an
open polydisc containing $\tau_{\infty}\in\mathfrak{D}_{f}.$ Let
$\mathcal{D}_{\tau_{\infty}}\subset\mathcal{U}_{\tau_{\infty}}\subset
\overline{\mathfrak{M}_{L}(\text{M})}$ be an open disk containing the point
$\tau_{\infty}\in\mathfrak{D}_{f}.$ Then the monodromy operator of the family
of polarized CY manifolds over $\mathcal{D}_{\tau_{\infty}}^{\ast}%
=\mathcal{D}_{\tau_{\infty}}-\tau_{\infty}$ is non-trivial group.
\end{proposition}

\textbf{Proof: }The subgroup $\Gamma$ of the mapping class group
$\Gamma\left(  \text{M}\right)  $ acts on the completion $\widetilde
{\mathcal{T}\left(  \text{M}\right)  }$ of the Teichm\"{u}ller space
$\mathcal{T}\left(  \text{M}\right)  $ with respect to the Hodge metric. Then
each of the components $\mathcal{D}_{i}$ of $\sigma^{-1}\left(  \Gamma
\left\backslash \widetilde{\mathcal{T}\left(  \text{M}\right)  }\right.
-\Gamma\left\backslash \mathcal{T}\left(  \text{M}\right)  \right.  \right)  $
is a fixed point set of some subgroup $G_{i}$ of $\Gamma.$ Let
\begin{equation}
\mathcal{X}_{\mathcal{D}_{\tau_{\infty}}}\rightarrow\mathcal{D}_{\tau_{\infty
}}\label{fam}%
\end{equation}
be a family of CY manifolds such that the $p(0)=\tau_{\infty}\in
\mathfrak{D}_{\infty}.$ Then the above arguments show that monodromy group of
the family is the stabilizer of the point $\tau_{\infty}\in\mathcal{D}_{i}.$
So Proposition \ref{NikA} is proved. $\blacksquare$

\begin{proposition}
\label{NikB}Let M be a CY manifold. Let $\mathcal{U}_{\tau_{\infty}}$ be an
open polydisc containing $\tau_{\infty}\in\mathfrak{D}_{f}.$ Let
$\mathcal{D}_{\tau_{\infty}}\subset\mathcal{U}_{\tau_{\infty}}\subset
\overline{\mathfrak{M}_{L}(\text{M})}$ be an open disk containing the point
$\tau_{\infty}\in\mathfrak{D}_{f}.$ Then then $\mathbf{L}^{2}$ metric on the
relative dualizing sheaf has growth \textbf{1. }$\left\vert \tau-\tau_{\infty
}\right\vert ^{1/k}$ if the monodoromy operator $T$ is of finite order on the
family $\left.  \mathcal{X}\right\vert _{\mathcal{D}_{\tau_{\infty}}%
-\tau_{\infty}}$, and \textbf{2. }$\left(  \log\left\vert \tau-\tau_{\infty
}\right\vert \right)  ^{k}$ if $T$ has an infinite order and $k>1.$
\end{proposition}

\textbf{Proof: }Let $\omega_{\tau}$ be a family of holomorphic $n-$forms over
the disk punctured disk $\mathcal{D}_{\tau_{\infty}}-\tau_{\infty}$. We need
to consider two cases:

\textbf{Case 1. }\textit{Suppose that the monodromy operator }$T$\textit{
acting on }$H^{n}\left(  \text{M}_{\tau},\mathbb{Z}\right)  $\textit{ is of
infinite order.} Let $\left\{  \gamma_{i}\right\}  $ be a basis of
$H_{n}\left(  \text{M}_{\tau},\mathbb{Z}\right)  $ for $\tau\neq\tau_{\infty
}.$ Let us consider
\begin{equation}
\left(  ...,a_{j}(\tau)=%
{\displaystyle\int\limits_{\gamma_{i}}}
\omega_{\tau},...\right)  .\label{A}%
\end{equation}
Then the components of $\left(  \ref{A}\right)  $ are solutions of a ordinary
differential equation with regular singular points. Let $\left\{  \gamma
_{i}^{\ast}\right\}  $ be the Dirac dual basis of $H^{n}\left(  \text{M}%
_{t},\mathbb{Z}\right)  ,$ i.e. $\left\langle \gamma_{i}^{\ast},\gamma
_{j}\right\rangle =\delta_{ij}.$ Then $\omega_{\tau}=%
{\displaystyle\sum\limits_{j=1}^{b_{n}}}
a_{j}(\tau)\gamma_{j}^{\ast}.$ It follows from Poincare duality that the
$\mathbf{L}^{2}$ metric is given by
\[
\left(  -1\right)  ^{\frac{n(n-1)}{2}}\left(  \frac{\sqrt{-1}}{2}\right)  ^{n}%
{\displaystyle\int\limits_{\text{M}}}
\omega_{\tau}\wedge\overline{\omega_{\tau}}=
\]%
\begin{equation}
\left(  ...,a_{j}(\tau)=%
{\displaystyle\int\limits_{\gamma_{i}}}
\omega_{\tau},...\right)  \left(  \left\langle \gamma_{i},\gamma
_{j}\right\rangle \right)  \left(  ...,a_{j}(\tau)=%
{\displaystyle\int\limits_{\gamma_{i}}}
\overline{\omega_{\tau}},...\right)  ^{t}.\label{mon1}%
\end{equation}

The assumption that the monodromy operator is of infinite order implies that
if $\gamma_{0}$ and $\gamma_{k}$ are cycles such that $T(\gamma_{0}%
)=\gamma_{0}$ and $\left(  T^{m}-id\right)  ^{k}\gamma_{k}=\gamma_{0}$ we get
that%
\begin{equation}
a_{k}(\tau)=%
{\displaystyle\int\limits_{\gamma_{k}}}
\omega_{\tau}\backsim\left(  \log\left(  \tau-\tau_{\infty}\right)  \right)
^{k}. \label{mon}%
\end{equation}
From the assumption that the monodromy operator is of infinite order, and from
$\left(  \ref{mon1}\right)  $ and $\left(  \ref{Disc3}\right)  $ we can
conclude that
\begin{equation}
\underset{\tau\rightarrow\tau_{\infty}}{\lim}\left\vert \frac{\left\Vert
\omega_{\tau}\right\Vert ^{2}}{\log|\tau-\tau_{\infty}|}\right\vert =c<\infty.
\label{Disc2}%
\end{equation}
So we proved Proposition \ref{NikB} in the case the monodromy is of infinite
oreder. So we need to consider the second case:

\textbf{Case 2. }\textit{Suppose that the monodromy operator }$T$\textit{
acting on the middle cohomology }$H^{n}\left(  \text{M}_{\tau}\text{,}%
\mathbb{Z}\right)  $\textit{ for }$\tau\in D_{\tau_{\infty}}-\tau_{\infty}%
$\textit{ has an order }$k>1.$ Let $\gamma_{0}$ and $\gamma_{1}\in
H_{n}\left(  \text{M,}\mathbb{Z}\right)  ,$ $\left\langle \gamma_{0}%
,\gamma_{k}\right\rangle \neq0$, $\gamma_{0}$ is an invariant vanishing cycle
of $T,$ i.e. $T(\gamma_{0})=\gamma_{0}$. Then according to \cite{LTY} we have
$T(\gamma_{k})\neq\gamma_{k}$ and $T^{k}(\gamma_{k})=\gamma_{k}.$ Then
according to the general theory of the monodromy we have that:%
\begin{equation}
\underset{\tau\rightarrow\tau_{\infty}}{\lim}%
{\displaystyle\int\limits_{\gamma_{0}}}
\omega_{\tau}=0\text{ and }\underset{\tau\rightarrow\tau_{\infty}}{\lim}\frac{%
{\displaystyle\int\limits_{\gamma_{1}}}
\omega_{\tau}}{\left(  \tau-\tau_{\infty}\right)  ^{1/k}}=c_{1}.\label{vc}%
\end{equation}
Formulas $\left(  \ref{mon}\right)  $ and $\left(  \ref{vc}\right)  $ imply
that
\begin{equation}
\underset{\tau\rightarrow\tau_{\infty}}{\lim}\frac{\left\Vert \omega_{\tau
}\right\Vert ^{2}}{\left\vert \tau-\tau_{\infty}\right\vert ^{1/k}}%
=c_{2}.\label{vc1}%
\end{equation}
So $\left(  \ref{Disc2}\right)  $ and $\left(  \ref{vc1}\right)  $ imply
Proposition \ref{NikB}. $\blacksquare$

\begin{proposition}
\label{NikC}Let M be a CY manifold. Let $\mathcal{U}_{\tau_{\infty}}$ be an
open polydisc containing $\tau_{\infty}\in\mathfrak{D}_{\infty}.$ Let
$\mathcal{D}_{\tau_{\infty}}\subset\mathcal{U}_{\tau_{\infty}}$ be an open
disk containing the point $\tau_{\infty}\in\mathfrak{D}_{\infty}.$ Then
$\overline{\eta^{N}}(\tau_{\infty})=0.$
\end{proposition}

\textbf{Proof: }Let $\mathcal{X}_{D}\rightarrow D$ be a family of CY manifolds
such that the $p(0)=\tau_{\infty}\in\mathfrak{D}_{\infty}.$ We need to
consider two cases:

\textbf{Case 1. }Suppose that the monodromy operator is of infinite order. Let
$\eta^{N}$ be the section constructed in Theorem \ref{Nik2} and $\left.
\eta^{N}\right\vert _{\mathcal{U}_{\tau_{\infty}}}=f_{\mathcal{U}%
_{\tau_{\infty}}}(\tau).$ From the explicit expression for the $\mathbf{L}%
^{2}$ norm of the section $\overline{\eta^{N}}$ we see that
\begin{equation}
\left.  \left\Vert \overline{\eta^{N}(\tau)}\right\Vert _{\mathbf{L}^{2}}%
^{2}\right\vert _{\mathcal{U}_{\tau}}=\frac{|f_{\mathcal{U}_{\tau_{\infty}}%
}(\tau)|^{2}}{\left\Vert \omega_{\tau}\right\Vert _{\mathbf{L}^{2}}^{2}}\leq
c_{0}, \label{Disc3}%
\end{equation}
where $\omega_{\tau}$ is a family of holomorphic $n-$forms defined over the
family of polarized CY manifolds obtained from the restriction of the versal
family $\mathcal{X}\rightarrow\mathfrak{M}_{L}\left(  \text{M}\right)  $ on
$\mathcal{U}_{\tau}$ and $\left.  \overline{\eta^{N}}\right\vert
_{\mathcal{U}_{\tau}}=f_{\mathcal{U}_{\tau}}(\tau).$ Thus $\left(
\ref{Disc3}\right)  $ and $\left(  \ref{Disc2}\right)  $ imply that%
\begin{equation}
\underset{\tau\rightarrow\tau_{\infty}}{\lim}f_{\mathcal{U}_{\tau}}(\tau)=0.
\label{Disc1}%
\end{equation}

\textbf{Case 2. }Suppose that the monodromy operator is of finite order. Let
$\eta^{N}$ be the section constructed in Theorem \ref{Nik2} and $\left.
\eta^{N}\right\vert _{\mathcal{U}_{\tau_{\infty}}}=f_{\mathcal{U}%
_{\tau_{\infty}}}(\tau).$ Then Formulas $\left(  \ref{Disc3}\right)  ,$
$\left(  \ref{vc1}\right)  ,$ and $\left\Vert \eta^{N}\right\Vert
_{\mathbf{L}^{2}}^{2}<\infty$ imply that $\underset{\tau\rightarrow
\tau_{\infty}}{\lim}f_{\mathcal{U}_{\tau_{\infty}}}(\tau)=0.$ Proposition
\ref{NikB} is proved. $\blacksquare$\ \

Propositions\ \ref{NikA} and \ref{NikB} imply Lemma \ref{NIK-}. $\blacksquare$
Lemma\ \ref{NIK-} implies Theorem \ref{Nik3}. $\blacksquare$

\begin{problem}
We constructed in Theorem \ref{Nik3} a holomorphic section $\overline{\eta
}^{N}$ of the prolonged line bundle $\overline{\omega_{\mathcal{X}%
/\mathfrak{M}_{L}(\text{M})}^{\otimes N}}$ for any dimension to some
compactification $\overline{\mathfrak{M}_{L}(\text{M})}$ of $\mathfrak{M}%
_{L}($M$)$\ such that $\overline{\mathfrak{M}_{L}(\text{M})}-\mathfrak{M}%
_{L}($M$)$ is a divisor with a normal crossings. Is it true that the function
defined by the regularized determinants $\det(\Delta_{0,1})$ of the Laplacians
of the CY metrics with imaginary parts cohomological to the polarization class
$L$ is bounded and $\left\Vert \eta^{\otimes N}\right\Vert _{\mathbf{L}^{2}%
}^{2}=\left(  \det(\Delta_{0,1})\right)  ^{N}$ and $N=\#\Gamma/\left[
\Gamma,\Gamma\right]  .$
\end{problem}

\subsection{The Analogue of Cusp Automorphic Forms}

\begin{definition}
\label{nef}Let M be a projective variety. A line bundle $\mathcal{L}$ on M
will be called big and nef if the Chern class $c_{1}\left(  \mathcal{L}%
\right)  $ of $\mathcal{L}$ satisfies the following conditions:
\end{definition}

\begin{enumerate}
\item \textit{For every irreducible curve }$C$\textit{ on M we have}
\begin{equation}%
{\displaystyle\int\limits_{C}}
c_{1}\left(  \mathcal{L}\right)  \geq0. \label{bn1}%
\end{equation}

\item \textit{The following inequality holds}%
\begin{equation}%
{\displaystyle\int\limits_{\text{M}}}
\wedge^{n}c_{1}\left(  \mathcal{L}\right)  >0. \label{bn2}%
\end{equation}

\end{enumerate}

\begin{definition}
\label{sa}A line bundle $\mathcal{L}$ on a complete scheme M is semi ample if
$\mathcal{L}^{\otimes m}$ is globally generated for some integer $m>0.$
\end{definition}

\begin{theorem}
\label{BB}Let M be a CY manifold. Let $\mathfrak{M}_{L}($M) be the moduli
space of polarized CY manifolds such that $\tau_{0}\in$ $\mathfrak{M}_{L}($M)
corresponds to M. Let $\omega_{\mathcal{X}/\mathfrak{M}_{L}(\text{M})}$ be the
relative dualizing sheaf of the family $\mathcal{X}\rightarrow\mathfrak{M}%
_{L}($M$).$ Let%
\[
H_{\infty}^{0}\left(  \mathfrak{M}_{L}(\text{M}),\omega_{\mathcal{X}%
/\mathfrak{M}_{L}(\text{M})}^{\otimes N}\right)  :=\left\{  s\in H_{\infty
}^{0}\left(  \mathfrak{M}_{L}(\text{M}),\omega_{\mathcal{X}/\mathfrak{M}%
_{L}(\text{M})}^{\otimes N}\right)  \left\vert \left\Vert s\right\Vert
_{\mathbf{L}^{2}}^{2}<\infty\right.  \right\}
\]
Then the space $H_{\infty}^{0}\left(  \mathfrak{M}_{L}(\text{M}),\omega
_{\mathcal{X}/\mathfrak{M}_{L}(\text{M})}^{\otimes N}\right)  $ is finite
dimensional for any $N$ and for large enough $N$ the linear system $H_{\infty
}^{0}\left(  \mathfrak{M}_{L}(\text{M}),\omega_{\mathcal{X}/\mathfrak{M}%
_{L}(\text{M})}^{\otimes N}\right)  $ defines a holomorphic inclusion
\[
\phi_{\infty}:\mathfrak{M}_{L}(\text{M})\subset\mathbb{CP}^{N_{0}}%
\]
such that $\phi_{\infty}\left(  \mathfrak{M}_{L}(\text{M})\right)  $ is a
Zariski open set.
\end{theorem}

\textbf{Proof: }Let $\overline{\mathfrak{M}_{L}(\text{M})}$ be any
compactification of $\mathfrak{M}_{L}($M$)$ such that
\[
\mathfrak{D}_{\infty}:=\overline{\mathfrak{M}_{L}(\text{M})}-\mathfrak{M}%
_{L}(\text{M})
\]
is a divisor of normal crossing. Let $\overline{\omega_{\mathcal{X}%
/\mathfrak{M}_{L}(\text{M})}}$ be the line bundle on $\overline{\mathfrak{M}%
_{L}(\text{M})}$ such that
\[
\overline{\omega_{\mathcal{X}/\mathfrak{M}_{L}(\text{M})}}\left\vert
_{\mathfrak{M}_{L}(\text{M})}\right.  =\omega_{\mathcal{X}/\mathfrak{M}%
_{L}(\text{M})}%
\]
and for any $\tau_{\infty}\in\mathfrak{D}_{\infty}$ and $U_{\tau_{\infty}}$
polydisk containing $\tau_{\infty}$ we have%
\[
\overline{\omega_{\mathcal{X}/\mathfrak{M}_{L}(\text{M})}}\left\vert
_{U_{\tau_{\infty}}}\right.  :=\left\{  s\in\Gamma\left(  U_{\tau_{\infty}%
}-U_{\tau_{\infty}}\cap\mathfrak{D}_{\infty},\omega_{\mathcal{X}%
/\mathfrak{M}_{L}(\text{M})}\right)  \left\vert \left\Vert s\right\Vert
_{\mathbf{L}^{2}}^{2}<\infty\right.  \right\}  .
\]
We will need the following Lemma:

\begin{lemma}
\label{BB1}The line bundle $\overline{\omega_{\mathcal{X}/\mathfrak{M}%
_{L}(\text{M})}}$ on $\overline{\mathfrak{M}_{L}(\text{M})}$ is a big, nef and semi-ample.
\end{lemma}

\textbf{Proof: }According to the Definition \ref{nef} we need to show that the
following Proposition:

\begin{proposition}
\label{BB1a}Let $C$ be any irreducible curve in $\overline{\mathfrak{M}%
_{L}(\text{M})}.$ Then we have:%
\[%
{\displaystyle\int\limits_{C}}
c_{1}\left(  \overline{\omega_{\mathcal{X}/\mathfrak{M}_{L}(\text{M})}%
},h\right)  \geq0.
\]

\end{proposition}

\textbf{Proof: }We proved in \cite{To89} that the Chern form of $\omega
_{\mathcal{X}/\mathfrak{M}_{L}(\text{M})}$ is the imaginary part of
Weil-Petersson metric of $\mathfrak{M}_{L}($M$).$ In \cite{To04} we proved
that the Weil-Petersson metric is a good metric in the sense of Mumford. See
\cite{Mum}. According to \cite{Mum} the Chern forms $c_{k}\left(  E,g)\text{
}\right)  $of a good metric $g$ on a vector bundle $E$ define currents which
are elements of $H^{2k}\left(  \mathfrak{M}_{L}(\text{M}),\mathbb{Z}\right)
.$ Let us denote by $h_{\mathbf{L}^{2}}$ the $\mathbf{L}^{2}$ metric on
$\omega_{\mathcal{X}/\mathfrak{M}_{L}(\text{M})}.$ Since the Chern form
$c_{1}\left(  \overline{\omega_{\mathcal{X}/\mathfrak{M}_{L}(\text{M})}%
},h_{\mathbf{L}^{2}}\right)  $ of the $\mathbf{L}^{2}$ metric $h_{\mathbf{L}%
^{2}}$ on $\omega_{\mathcal{X}/\mathfrak{M}_{L}(\text{M})}$ is the imaginary
part of the Weil-Petersson metric and since $h_{\mathbf{L}^{2}}$ is a good
metric, then the current $c_{1}\left(  \overline{\omega_{\mathcal{X}%
/\mathfrak{M}_{L}(\text{M})}},h_{\mathbf{L}^{2}}\right)  $ is positive. This
fact implies that for any irreducible curve $C$ on $\overline{\mathfrak{M}%
_{L}(\text{M})}$ we have
\[%
{\displaystyle\int\limits_{C}}
c_{1}\left(  \overline{\omega_{\mathcal{X}/\mathfrak{M}_{L}(\text{M})}%
},h_{\mathbf{L}^{2}}\right)  \geq0.
\]
The last inequality implies Proposition \ref{BB1a}. $\blacksquare$

Since $h_{\mathbf{L}^{2}}=\mathbf{L}^{2}$ is a good metric on $\omega
_{\mathcal{X}/\mathfrak{M}_{L}(\text{M})}$\ and $dd^{c}\log h_{\mathbf{L}^{2}%
}$ is the imaginary part of the Weil-Petersson metric, then we have%
\[%
{\displaystyle\int\limits_{\overline{\mathfrak{M}_{L}(\text{M})}}}
\wedge^{h^{n-1,1}}c_{1}\left(  \overline{\omega_{\mathcal{X}/\mathfrak{M}%
_{L}(\text{M})}},h_{\mathbf{L}^{2}}\right)  >0.
\]
So $\overline{\omega_{\mathcal{X}/\mathfrak{M}_{L}(\text{M})}}$ is a big and
nef. According to Theorem \ref{sec}, there exists a non-zero holomorphic
section $\eta^{\otimes N}.$ The definition of $\overline{\omega_{\mathcal{X}%
/\mathfrak{M}_{L}(\text{M})}}$ implies that $\eta^{\otimes N}$ generates
$\overline{\omega_{\mathcal{X}/\mathfrak{M}_{L}(\text{M})}^{\otimes N}}.$ So
$\omega_{\mathcal{X}/\mathfrak{M}_{L}(\text{M})}$ is semiample nef and big
line bundle. Theorem Lemma \ref{BB1} is proved. $\blacksquare$

Theorem \textbf{2.1.27 }proved on page 129 of \cite{L} implies Theorem
\ref{BB}. Before formulating Theorem \textbf{2.1.27 }we will introduce some
notions. Let L be some semi-ample line bundle. We will denote by $M\left(
X,L\right)  $ the semi group
\[
M\left(  X,\mathcal{L}\right)  =\left\{  \left.  m\in\mathbb{N}\right\vert
\mathcal{L}^{\otimes m}\text{ has no fixed points}\right\}  .
\]
$f(\mathcal{L})$ will be the exponent of $M\left(  X,\mathcal{L}\right)  ,$
i.e. the largest natural number such that every element of $M\left(
X,\mathcal{L}\right)  $ is a multiple of $g(\mathcal{L}).$ Given $m\in
M\left(  X,\mathcal{L}\right)  $ we will denote by $Y_{m}=\phi_{m}(X)$ the
image of $X$ by the holomorphic map%
\[
\phi_{m}:X\rightarrow\phi_{m}(X)=Y_{m}\subset\mathbb{P}\left(  H^{0}\left(
X,\mathcal{L}^{\otimes m}\right)  \right)  ,
\]
determined by the linear system $\left\vert \mathcal{L}^{\otimes m}\right\vert
.$

\textbf{Theorem 2.1.27 (semi-ample fibrations). }\textit{Let X be a normal
projective space and let }$L$\textit{ be a semi-ample line bundle on }%
$X.$\textit{ Then there is a fibre space }$\phi:X\rightarrow Y$\textit{ having
the property that for any sufficiently large }$m\in M(X,L),$\textit{ }%
$Y_{m}=Y$\textit{ and }$\phi_{m}=\phi$\textit{. Moreover there is an ample
line bundle }$A$\textit{ on }$Y$\textit{ such that} $\mathcal{L}=\phi^{\ast
}(A).$ Theorem \ref{BB} is proved. $\blacksquare$

\subsection{Extension Properties}

\begin{theorem}
\label{B}Let $\phi:\left(  D^{\ast}\right)  ^{k}\times D^{m}\rightarrow
\mathfrak{M}_{L}($M$)$ be a holomoprphic map. Then $\phi$ can be extended to a
holomorphic map $\overline{\phi}:D^{k}\times D^{m}\rightarrow\widetilde
{\mathfrak{M}_{L}(\text{M})},$ where $\widetilde{\mathfrak{M}_{L}(\text{M})}$
is the Baily-Borel compactification of $\mathfrak{M}_{L}($M$).$
\end{theorem}

\textbf{Proof: }The proof of Theorem \ref{B} is based on the following
generalization of Schwarz Lemma due to S. -T. Yau;

\begin{theorem}
\label{yausch}Let M be a complete K\"{a}hler manifold with a Ricci curvature
bounded from below by a constant, and N be another Hermitian manifold with
holomorphic bisectional curvature bounded from above by a negative constant.
Then any holomorphic mapping from M into N decreases distances up to a
constant depending only on the curvature of M and N. See \cite{yau}.
\end{theorem}

Theorem \ref{B} follows from the following Lemma:

\begin{lemma}
\label{B1}Let $\phi:\left(  D^{\ast}\right)  ^{k}\times D^{m}\rightarrow
\mathfrak{M}_{L}($M$)$ be a holomoprphic map. Let $(z_{k},w_{k})\in\left(
D^{\ast}\right)  ^{k}\times D^{m}$ such
\[
\underset{k\rightarrow\infty}{\lim}(z_{k},w_{k})=\left(  0,w_{0}\right)  \in
D^{k}\times D^{m}.
\]
Then $\underset{k\rightarrow\infty}{\lim}\phi(z_{k},w_{k})$ exists and%
\[
\underset{k\rightarrow\infty}{\lim}\phi(z_{k},w_{k})\in\widetilde
{\mathfrak{M}_{L}(\text{M})}%
\]
$\widetilde{\mathfrak{M}_{L}(\text{M})},$ where $\widetilde{\mathfrak{M}%
_{L}(\text{M})}$ is the Baily-Borel compactification of $\mathfrak{M}_{L}%
($M$).$
\end{lemma}

\textbf{Proof: }Let us consider on $\left(  D^{\ast}\right)  ^{k}\times D^{m}$
the Poincare metric:%
\[%
{\displaystyle\sum\limits_{i=1}^{k}}
\frac{dz^{i}\otimes\overline{dz^{i}}}{\left\vert z^{i}\right\vert ^{2}\left(
\log\left\vert z^{i}\right\vert \right)  ^{2}}+%
{\displaystyle\sum\limits_{j=1}^{m}}
\frac{dw^{j}\otimes\overline{dw^{j}}}{\left(  1-\left\vert w^{j}\right\vert
^{2}\right)  ^{2}}.
\]
Then to apply Theorem \ref{yausch} we replace M with $\left(  D^{\ast}\right)
^{k}\times D^{m}$ with Poincare metric and replace N with $\widetilde
{\mathfrak{M}_{L}(\text{M})}$ with Hodge metric. Then Lemma \ref{B1} follows
from Theorem \ref{yausch} directly since the holomorphic map $\phi$ is a
distance decreasing up to a constant, with respect to the Hodge metric on
$\mathfrak{M}_{L}($M$)$. Theorem \ref{NC} implies that the Hodge metric
satisfies the conditions of Theorem \ref{yausch}. See \cite{LSTY}.
$\blacksquare$

Lemma \ref{B1} implies that the map $\phi:\left(  D^{\ast}\right)  ^{k}\times
D^{m}\rightarrow\mathfrak{M}_{L}($M$)$ can be extended to a continuous map:%
\[
\overline{\phi}:D^{k}\times D^{m}\rightarrow\widetilde{\mathfrak{M}%
_{L}(\text{M})}.
\]
Then according to Riemann extension Theorem (See Theorem \textbf{44.42} on p.
420 in \cite{Ab}) $\overline{\phi}$ is a holomorphic map. Theorem \ref{B} is
proved. $\blacksquare$

\begin{corollary}
\label{BA}Let $\mathcal{Z}$ be a quasi-projective variety. Suppose that%
\[
\phi:\mathcal{Z}\rightarrow\mathfrak{M}_{L}(\text{M})
\]
be a non trivial holomorphic map. Let $\overline{\mathcal{Z}}$ be a projective
manifold such that $\overline{\mathcal{Z}}-\mathcal{Z}$ is a divisor of normal
crossings. Then $\phi$ can be extended to a holomorphic map%
\[
\overline{\phi}:\overline{\mathcal{Z}}\rightarrow\widetilde{\mathfrak{M}%
_{L}(\text{M})}.
\]

\end{corollary}

\begin{remark}
Theorem \ref{B} is a generalization of a Theorem of A. Borel proved for
locally symmetric Hermitian spaces in \cite{ab}. Theorem \ref{B} implies that
$\widetilde{\mathfrak{M}_{L}(\text{M})}$ is a minimal model among all possible
compactifications of the moduli space $\mathfrak{M}_{L}($M$)$ with a boundary
divisors with normal crossings.
\end{remark}

\section{CY\ Threefolds}

\subsection{Invariants of the Short Term Asymptotic Expansion of the Heat
Kernel}

\begin{theorem}
\label{And1}Suppose that M is a three dimensional CY manifold and g is a CY
metric. Then the coefficients a$_{2k}$ for $k=3,2,1$ and $0$ in the expression
$\left(  \ref{ass0}\right)  $ are constants which depend only on the CY
manifolds and the fixed class of cohomology of the CY metric.
\end{theorem}

\textbf{Proof:} We know that the Heat kernel has the following asymptotic
expansion:
\begin{equation}
Tr(\exp(-t\triangle_{0}))=\frac{a_{-n}(g)}{t^{n}}+\frac{a_{-n+1}(g)}{t^{n-1}%
}+\frac{a_{-n+2}(g)}{t^{n-2}}+...+a_{0}(g)+h(t,\tau,\overline{\tau}).
\label{ass0}%
\end{equation}
(See \cite{Roe}.) We will apply $\left(  \ref{ass0}\right)  $ for three
dimensional CY manifolds. In \cite{Gil} on page 118 one can find the following
formulas for $a_{-3}(g),a_{-2}(g)$ and $a_{-1}(g):$%
\[
\alpha_{-3}(g)=\frac{vol(g)}{4\pi},\text{ }a_{-2}(g)=\frac{-\int
\limits_{\text{M}}k(g)vol(g)}{24\pi}%
\]
and
\begin{equation}
a_{-1}(g)=\frac{-12\left(  \int\limits_{\text{M}}\Delta_{g}%
(k(g))vol(g)\right)  +5\left\Vert Ric(g)\right\Vert ^{2}-2\left\Vert
R(g)\right\Vert ^{2}}{1440\pi} \label{ass1}%
\end{equation}
where $k(g)$ is the scalar curvature of the metric $g,$ $\left\Vert
Ric(g)\right\Vert ^{2}$ is the $L^{2}$ norm of the Ricci tensor of $g$ and
$\left\Vert R(g)\right\Vert ^{2}$ is the $L^{2}$ norm of the curvature of the
metric $g.$ Using the fact that $g$ is a Calabi-Yau metric, i.e.
$Ric(g)=k(g)=0,$ we obtain:
\begin{equation}
\alpha_{-n}(g)=\frac{vol(g)}{4\pi},\text{ }a_{-1}(g)=0\text{ and }%
a_{-2}(g)=\frac{-\left\Vert R(g)\right\Vert ^{2}}{720\pi}. \label{ass2}%
\end{equation}
In \cite{Ca} Calabi proved on page 264 the following Proposition:

\begin{proposition}
\label{Ca}The following formula holds on a complex K\"{a}hler manifold M with
a fixed cohomology class $L$ of the imaginary part of a K\"{a}hler metric:
\begin{equation}
2\left\|  Ric(g)\right\|  ^{2}-\left\|  R(g)\right\|  ^{2}-\int
\limits_{\text{M}}k(g)^{2}vol(g)=-\int\limits_{\text{M}}c_{2}(\text{M}%
)\wedge\omega_{g}^{n-2}, \label{ass3}%
\end{equation}
where $c_{2}($M$)$ is the second Chern class of M.
\end{proposition}

Applying formula $\left(  \ref{ass3}\right)  $ to a CY metric, we obtain that
on a three dimensional CY manifold with a Calabi Yau metric $g$ we have:
\begin{equation}
a_{-n}(g)=\frac{1}{4\pi}\int\limits_{\text{M}}L^{n},\text{ }a_{n-1}(g)=0\text{
and }a_{n-2}(g)=-\frac{1}{720\pi}\int\limits_{\text{M}}c_{2}(\text{M})\wedge
L^{n-2}. \label{ass4}%
\end{equation}
Theorem \ref{And1} follows directly from $\left(  \ref{ass4}\right)  $ since
$\left(  \ref{ass4}\right)  $ implies that $a_{-1}(g),$ $a_{-2}(g)$ and
$a_{-3}(g)$ are topological invariants$.$ We need to prove that $a_{0}(g)$ is
a constant in order to deduce Theorem \ref{And1}.

\begin{lemma}
\label{Ca1}Let%
\[
Tr(\exp(-t\Delta_{0})=\frac{a_{-n}}{t^{n}}+...+\frac{a_{-1}}{t}+a_{0}+O(t)
\]
be the asymptotic expansion of $Tr(\exp(-t\Delta_{0})$ with respect to a CY
metric with a fixed class of cohomology of its imaginary part. Then the real
coefficient $a_{0}$ is a constant, i.e.
\[
\frac{\partial}{\partial\tau}a_{0}(\tau,\overline{\tau})=\frac{\overline
{\partial}}{\overline{\partial\tau}}a_{0}(\tau,\overline{\tau})=0.
\]

\end{lemma}

\textbf{Proof: }According to \cite{BGV} the following equality is true:
\begin{equation}
\zeta_{0,\tau}(0)=a_{0}(\tau), \label{ass5a}%
\end{equation}
where $a_{0}$ is a real valued function on the moduli space of polarized CY
manifolds. If we prove that:
\begin{equation}
\left.  \frac{\partial}{\partial\tau}\left(  \zeta_{0,\tau}(s)\right)
\right\vert _{s=0}=0, \label{ass5}%
\end{equation}
then Lemma \ref{Ca1} will follow directly from $\left(  \ref{ass5a}\right)  $.
In the paper \cite{BT} on page 85 the following formula was proved:
\[
\frac{\overline{\partial}}{\overline{\partial\tau^{i}}}\left(  \zeta_{0,\tau
}^{"}(s)\right)  =
\]%
\[
\frac{1}{\Gamma(s)}\int\limits_{0}^{\infty}Tr\left(  \exp\left(
-t\triangle_{\tau,0}^{^{"}}\right)  \circ\triangle_{\tau}^{^{"}}\circ\left(
\partial_{\tau}\right)  ^{-1}\circ F^{^{\prime}}(1,\overline{\frac{\partial
}{\partial\tau^{i}}\phi(\tau)})\circ\overline{\partial_{\tau}}\right)
t^{s}dt.
\]
repeating word by word the proof the above formula from \cite{BT} we get:%
\[
\frac{\partial}{\partial\tau^{i}}\left(  \zeta_{0,\tau}^{"}(s)\right)  =
\]%
\begin{equation}
\frac{1}{\Gamma(s)}\int\limits_{0}^{\infty}Tr\left(  \frac{d}{dt}\left(
\exp\left(  -t\triangle_{\tau,0}^{^{"}}\right)  \right)  \circ\left(
\overline{\partial}_{\tau}\right)  ^{-1}\circ F^{^{\prime}}(q,\frac{\partial
}{\partial\tau^{i}}\phi(\tau))\circ\partial_{\tau}\right)  t^{s}dt.
\label{ass6}%
\end{equation}
By integrating by parts the expressions in $\left(  \ref{ass6}\right)  ,$ we
obtain:
\[
\frac{\partial}{\partial\tau^{i}}\left(  \zeta_{0,\tau}^{"}(s)\right)  =
\]%
\begin{equation}
\frac{s}{\Gamma(s)}\int\limits_{0}^{\infty}Tr\left(  \exp\left(
-t\triangle_{\tau,0}^{^{"}}\right)  \circ\left(  \overline{\partial}_{\tau
}\right)  ^{-1}\circ F^{^{\prime}}(q,\frac{\partial}{\partial\tau^{i}}%
\phi(\tau))\circ\partial_{\tau}\right)  t^{s-1}dt. \label{ass7}%
\end{equation}
We can rewrite the integral in the right hand side of $\left(  \ref{ass7}%
\right)  $ as follows:
\[
\frac{s}{\Gamma(s)}\int\limits_{0}^{\infty}Tr\left(  \exp\left(
-t\triangle_{\tau,0}^{^{"}}\right)  \circ\left(  \overline{\partial}_{\tau
}\right)  ^{-1}\circ F^{^{\prime}}(q,\frac{\partial}{\partial\tau^{i}}%
\phi(\tau))\circ\partial_{\tau}\right)  t^{s-1}dt=
\]%
\[
\frac{s}{\Gamma(s)}\int\limits_{0}^{1}Tr\left(  \exp\left(  -t\triangle
_{\tau,0}^{^{"}}\right)  \circ\left(  \overline{\partial}_{\tau}\right)
^{-1}\circ F^{^{\prime}}(q,\frac{\partial}{\partial\tau^{i}}\phi(\tau
))\circ\partial_{\tau}\right)  t^{s-1}dt+
\]%
\begin{equation}
\frac{s}{\Gamma(s)}\int\limits_{1}^{\infty}Tr\left(  \exp\left(
-t\triangle_{\tau,0}^{^{"}}\right)  \circ\left(  \overline{\partial}_{\tau
}\right)  ^{-1}\circ F^{^{\prime}}(q,\frac{\partial}{\partial\tau^{i}}%
\phi(\tau))\circ\partial_{\tau}\right)  t^{s-1}dt. \label{ass7a}%
\end{equation}
From the short term asymptotic expansion
\[
Tr\left(  \exp\left(  -t\triangle_{\tau,0}^{^{"}}\right)  \circ\left(
\overline{\partial}_{\tau}\right)  ^{-1}\circ F^{^{\prime}}(q,\frac{\partial
}{\partial\tau^{i}}\phi(\tau))\circ\partial_{\tau}\right)  =
\]%
\begin{equation}
\frac{c_{-k}(\tau)}{t^{k}}+...+\frac{c_{-1}(\tau)}{t}+c_{0}(\tau)+O(t)
\label{ass7b}%
\end{equation}
we obtain that
\[
\frac{s}{\Gamma(s)}\int\limits_{0}^{1}Tr\left(  \exp\left(  -t\triangle
_{\tau,0}^{^{"}}\right)  \circ\left(  \overline{\partial}_{\tau}\right)
^{-1}\circ F^{^{\prime}}(q,\frac{\partial}{\partial\tau^{i}}\phi(\tau
))\circ\partial_{\tau}\right)  t^{s-1}dt=
\]%
\begin{equation}
\frac{s}{\Gamma(s)}\left(  \frac{c_{0}(\tau)}{s}+\gamma_{0}(\tau)+O(s)\right)
. \label{ass7c}%
\end{equation}
From $\left(  \ref{ass7c}\right)  $ and the fact that $\frac{s}{\Gamma
(s)}=s^{2}+O(s^{3})$ we get that
\[
\frac{\partial}{\partial\tau^{i}}\left(  \zeta_{0,\tau}^{"}(s)\right)
=s^{2}\left(  \frac{c_{0}(\tau)}{s}+\gamma_{0}(\tau)+O(s)\right)  .
\]
From the last formula we obtain that
\[
\left(  \frac{\partial}{\partial\tau^{i}}\left(  \zeta_{0,\tau}^{"}\right)
\right)  (0)=0.
\]
Lemma \ref{Ca1} is proved. $\blacksquare$

Lemma \ref{Ca1} implies Theorem \ref{And1}. $\blacksquare$

\subsection{The Regularized CY Determinants are Bounded}

\begin{theorem}
\label{LY}For CY threefolds the regularized determinants of the Laplacians
$\Delta_{q}$ of the Calabi Yau metrics $g(\tau,\overline{\tau})$ with a fixed
cohomology class $L$ for $0\leq q\leq n=\dim_{\mathbb{C}}$M are bounded as
functions on the moduli space, i.e. we have: $0\leq\det(\triangle_{\tau
,q})\leq C_{q}.$
\end{theorem}

\textbf{Proof: }We will outline the main ideas of the proof of Theorem
\ref{LY}. In order to prove Theorem \ref{LY} it is enough to bound
$\det(\triangle_{0}).$ The bound of $\det(\triangle_{0})$ is based on the
following expression for the zeta function of the Laplacian acting on functions:%

\[
\zeta_{0}(s)=\frac{1}{\Gamma(s)}\int\limits_{\text{0}}^{\infty}\left(
Tr(\exp(-t\triangle_{0}\right)  )t^{s-1}dt=b_{0}+b_{1}s+O(s^{2}).
\]
From the definition of $\det(\triangle_{0})$ it follows that
\begin{equation}
\det(\triangle_{0})=\exp(-b_{1}). \label{LY2}%
\end{equation}
So if $b_{1}$ is bounded from bellow, i.e.
\begin{equation}
c_{1}\leq b_{1} \label{LY3}%
\end{equation}
then Theorem \ref{LY} will be proved. The bound of $b_{1}$ is based on several facts.

\begin{enumerate}
\item The following explicit formula for $b_{1}$ is proved in \cite{ABKS}:
\begin{equation}
b_{1}=\gamma a_{0}+\sum_{k=1}^{3}\frac{a_{-k}}{k}+\psi_{1}+\psi_{2},
\label{LY1}%
\end{equation}
where $\gamma$ is the Euler constant, $\psi_{1}$ is given by the formula
\[
\psi_{1}(t,\tau,\overline{\tau})=\int\limits_{0}^{1}\left(  Tr(\exp
(-t\triangle_{0}))-\sum_{k=0}^{3}\frac{a_{-k}}{t^{k}}\right)  \frac{dt}{t}%
\]
and $\psi_{2}$ by
\[
\psi_{2}(t,\tau,\overline{\tau})=\int\limits_{1}^{\infty}Tr(\exp
(-t\triangle_{0}))\frac{dt}{t}.
\]

\item Theorem \ref{And1} implies that the expression: $\gamma a_{0}+\sum
_{k=1}^{3}\frac{a_{-k}}{k}$ in $\left(  \ref{LY1}\right)  $ is a constant.
Clearly $\psi_{2}(t,\tau,\overline{\tau})>0.$

\item The third fact is that $\psi_{1}(t,\tau,\overline{\tau})\geq c_{0},$
where $c_{0}$ is a constant$.$
\end{enumerate}

Combining all these facts we get $\left(  \ref{LY3}\right)  $. $\left(
\ref{LY3}\right)  $ and the explicit formula $\left(  \ref{LY1}\right)  $ will
imply that $0\leq\det(\Delta_{0})\leq C<\infty.$ So we need to prove the
following Lemma:

\begin{lemma}
\label{end}The following inequality holds: $\psi_{1}(t,\tau,\overline{\tau
})\geq c_{0}.$
\end{lemma}

\textbf{Proof: }Let
\begin{equation}
h(t,\tau,\overline{\tau})=Tr(\exp(-t\triangle_{0}))-\sum_{k=0}^{3}\frac
{a_{-k}}{t^{k}} \label{b100}%
\end{equation}
We also know that $h(t,\tau,\overline{\tau})=th_{1}(t,\tau,\overline{\tau}).$
According to Theorem \ref{And1}, the expression $%
{\displaystyle\sum\limits_{i=1}^{3}}
\frac{a_{-k}}{t^{k}}$ is a function which does not depend on $\tau$ and
$\overline{\tau}.$ We also know that $Tr(\exp(-t\triangle_{\tau,0}))$ is a
strictly positive function for $t>0$ which depends on $t,\tau$ and
$\overline{\tau}.$ Thus $\underset{\tau}{\inf}\left(  Tr(\exp(-t\triangle
_{\tau,1}))\right)  $ exists for $t>0$. Let
\[
\phi(t):=\underset{\tau}{\inf}\left(  Tr(\exp(-t\triangle_{\tau,1}))\right)
-\sum_{k=0}^{3}\frac{a_{-k}}{t^{k}}%
\]
for $\,0<t\leq1.$ Then
\[
\phi(t)=\underset{\tau}{\inf}\left(  Tr(\exp(-t\triangle_{\tau,1}))\right)
-\sum_{k=0}^{3}\frac{a_{-k}}{t^{k}}=
\]%
\[
\underset{\tau}{\inf}\left(  Tr(\exp(-t\triangle_{\tau,1}))-\sum_{k=0}%
^{3}\frac{a_{-k}}{t^{k}}\right)  =\underset{\tau}{\inf}\left(  h(t,\tau
,\overline{\tau})\right)  .
\]

\begin{proposition}
\label{end1}$\phi(t)$ is a continuous function for $0\leq t\leq1.$
\end{proposition}

\textbf{Proof: }Let $0\leq t_{0}\leq1$ be a fixed number. Let $\left\{
t_{n}\right\}  $ be a sequence of numbers such that $0<t_{n}<1$ such that
$\underset{n\rightarrow\infty}{\lim}t_{n}=t_{0}$ then we need to prove that
\begin{equation}
\underset{n\rightarrow\infty}{\lim}\phi(t_{n})=\underset{n\rightarrow\infty
}{\lim}\left(  \underset{\tau}{\inf}\left(  h(t,\tau,\overline{\tau})\right)
\right)  =\underset{\tau}{\inf}\left(  h(t,\tau,\overline{\tau})\right)
=\phi(t_{0}). \label{E1}%
\end{equation}
The definition of the function $h(t,\tau,\overline{\tau})$ show that
$h(t,\tau,\overline{\tau})$ is a continuous function such that\ $h(0,\tau
,\overline{\tau})=0.$ The definition of $\underset{\tau}{\inf}h(t,\tau
,\overline{\tau})$ implies that there exists a sequence $\left\{  \tau
_{n}\right\}  $ of points $\tau_{n}\in\mathfrak{M}_{L}($M$)$ such that
\begin{equation}
\underset{n\rightarrow\infty}{\lim}\left(  h(t,\tau,\overline{\tau})\right)
=\underset{\tau}{\inf}\left(  h(t,\tau,\overline{\tau})\right)  . \label{E2}%
\end{equation}
Since $h(t,\tau,\overline{\tau})$ is a continuous function we get that
\begin{equation}
\underset{n\rightarrow\infty}{\lim}\left(  \underset{k\rightarrow\infty}{\lim
}\left(  h(t_{k},\tau_{n},\overline{\tau_{n}})\right)  \right)  =\underset
{n\rightarrow\infty}{\lim}\left(  h(t_{0},\tau_{n},\overline{\tau_{n}%
})\right)  =\underset{\tau}{\inf}\left(  h(t,\tau,\overline{\tau})\right)  .
\label{E3}%
\end{equation}
On the other hand we have%
\begin{equation}
\underset{n\rightarrow\infty}{\lim}\left(  \underset{k\rightarrow\infty}{\lim
}\left(  h(t_{k},\tau_{n},\overline{\tau_{n}})\right)  \right)  =\underset
{k\rightarrow\infty}{\lim}\left(  \underset{n\rightarrow\infty}{\lim}\left(
h(t_{k},\tau_{n},\overline{\tau_{n}})\right)  \right)  =\underset
{k\rightarrow\infty}{\lim}\left(  \underset{\tau}{\inf}\left(  h(t_{k}%
,\tau,\overline{\tau})\right)  \right)  . \label{E4}%
\end{equation}
Combining $\left(  \ref{E3}\right)  $ and $\left(  \ref{E4}\right)  $ we get
that
\[
\underset{k\rightarrow\infty}{\lim}\left(  \underset{\tau}{\inf}\left(
h(t_{k},\tau,\overline{\tau})\right)  \right)  =\underset{k\rightarrow\infty
}{\lim}\phi(t_{k})=\underset{\tau}{\inf}\left(  h(t_{0},\tau,\overline{\tau
})\right)  =\phi(t_{0}).
\]
Proposition \ref{end1} is proved. $\blacksquare$

Since $\phi(t)$ is a continuous function on the closed interval $[0,1]$ it has
a minimum. Let $c_{0}=\underset{0\leq t\leq1}{\min}\phi(t).$ So we have
\[
\int\limits_{\text{0}}^{1}(h(t,\tau,\overline{\tau})-c_{0})dt\geq0.
\]
On the other hand, we have
\[
\int\limits_{\text{0}}^{1}(h(t,\tau,\overline{\tau})-c_{0})dt=\int
\limits_{\text{0}}^{1}h(t,\tau,\overline{\tau})dt-c_{0}=\psi_{1}%
(t,\tau,\overline{\tau})-c_{0}\geq0.
\]
So $\psi_{1}(t,\tau,\overline{\tau})\geq c_{0}.$ Lemma \ref{end} is proved.
$\blacksquare$

Theorem \ref{LY} is proved. $\blacksquare$

\subsection{The Existence of Global Section $\eta^{\otimes N}$ of the Relative
Dualizing Sheaf with Finite $\mathbf{L}^{2}$ Norm $\left(  \det(\Delta
_{0,1})\right)  ^{N}$}

\begin{theorem}
\label{Nik1a}Let M be a three dimensional CY manifold and let $N=\#\Gamma
/[\Gamma,\Gamma].$ Then there exists a section $\eta^{N}$ of the line bundle
$\overline{\left(  \omega_{\mathcal{X}/\mathfrak{M}_{L}(\text{M})}^{\otimes
N}\right)  }$ such that the norm of $\eta^{N}$ with respect to the $N$ tensor
power of the $L^{2}$ metric on $\omega_{\mathcal{X}/\mathfrak{M}_{L}%
(\text{M})}^{\otimes N}$ is given by: $\left\Vert \eta^{N}\right\Vert
_{\mathbf{L}^{2}}^{2}=\left(  \det(\Delta_{0,1})\right)  ^{N}.$ The zero set
of $\eta^{N}$ is a non zero effective divisor whose support contains or is
equal to the support of $\mathfrak{D}_{\infty}$, where $\mathfrak{D}_{\infty}$
is defined in Definition \ref{Disc}.
\end{theorem}

\textbf{Proof: }Theorem \ref{sec} implies that there exists a the section
$\eta^{\otimes N}$ of the relative dualizing line bundle $\omega
_{\mathcal{X}/\mathfrak{M}_{L}(\text{M})}^{\otimes N}$ on $\mathfrak{M}_{L}%
($M$)$ such that $\left\Vert \eta(\tau)\right\Vert _{\mathbf{L}^{2}}^{2}%
=\det\Delta_{(0,1)}(\tau)$ and $\left\Vert \eta(\tau)\right\Vert
_{\mathbf{L}^{2}}^{2}>0$ for $\tau\in\mathfrak{M}_{L}($M$)\mathbf{.}$ Now we
can apply Propositions \ref{NikA}, \ref{NikB} and \ref{NikC} to conclude the
zero set $\left(  \overline{\eta^{N}}\right)  _{0}$ is the effective divisor
$\mathfrak{D}_{\infty}$ defined by Definition \ref{Disc}. Theorem\ \ref{Nik1a}
is proved. $\blacksquare.$


\begin{thebibliography}{99}                                                                                               %


\bibitem {Ab}S. Abhyankar, "\textit{Local Analytic Geometry", }Academic Press,
New York, 1964.

\bibitem {ABKS}D. Abramovich, J. -F. Burnol, J. Kramer and C. Soul\'{e},
\textquotedblright\textit{Lectures on Arakelov Geometry\textquotedblright},
Cambridge Studies In Advanced Mathematics Volume \textbf{33}, Cambridge
University Press, 1992.

\bibitem {BCOV}M. Bershadsky, S. Cecotti, H. Ooguri and C. Vafa,
\textit{\textquotedblright Kodaira-Spencer Theory of Gravity and Exact Results
for Quantum String Amplitude\textquotedblright}, Comm. Math. Phys.
\textbf{165} (1994), 311-428.

\bibitem {BT}J. Bass, A. Todorov, "\textit{The Analogue of the Dedekind Eta
Function for CY Manifolds I", }Journal fur die reine und angewandte Mathematik
(Crelles Journal), v. 599, 61-96 (2006).

\bibitem {ab}A. Borel, "\textit{Some Metric Properties of Arithmetic Quotients
of Symmetric Spaces and Extension Theorem", }J. Diff. Geometry, \textbf{vol. 6
}(1972), 543-560.

\bibitem {Ca}E. Calabi, \textquotedblright\textit{Extremal K\"{a}hler
Metrics\textquotedblright, }Seminar on Differential Geometry, ed. S.-T. Yau,
Annals of Mathematical Studies, vol. \textbf{102}, Princeton University Press.

\bibitem {BGV}N. Berline, E. Getzler and \ M. Vergne,
\textit{\textquotedblright Heat Kernels and Dirac Operators\textquotedblright,
}Springer Verlag, 1991.

\bibitem {Bour}A. A. Kirillov and Cl. Delaroche, \textquotedblright\textit{Sur
les Relations entre l'espace dual d'un group et la struture de ses
sous-groupes ferm\'{e}s\textquotedblright\ (d'apr\`{e}s D. A. Kazhdan)},
Bourbaki No. 343, 1967/68.

\bibitem {D}S. Donaldson, \textquotedblright\textit{Infinite Determinants,
Stable Bundles and Curvature\textquotedblright, }Duke Mathematical Journal,
vol. \textbf{54, Number 1}(1987), 231-248.

\bibitem {DK}S. Donaldson and P. Kronheimer, \textquotedblright\textit{The
Geometry of Four Manifolds\textquotedblright,} Oxford Mat. Monog., Oxford
Science Publications, Oxford University Press, New York 1990.

\bibitem {Gil}P. Gilkey, \textquotedblright\textit{Invariance Theory, The Heat
Equation, And the Atiyah-Singer Index Theorem\textquotedblright,} Mathematics
Lecture Series vol. \textbf{11}, Publish or Parish, Inc. Wilmington, Delaware
(USA) 1984.

\bibitem {JT95}J. Jorgenson and A. Todorov, \textquotedblright\textit{A
Conjectural Analogue of Dedekind Eta Function for K3
Surfaces\textquotedblright, }\ Math. Research Lett. \textbf{2}(1995) 359-360.

\bibitem {JT98}J. Jorgenson and A. Todorov, \textquotedblright\textit{Analytic
Discriminant for Polarized Algebraic K3 Surfaces\textquotedblright,} Mirror
Symmetry III, ed. S-T. Yau and Phong, AMS, p. 211-261.

\bibitem {Mum}D. Mumford, \textquotedblright\textit{Hirzebruch's
Proportionality Principle in the Non-Compact Case\textquotedblright,} Inv.
Math. \textbf{42}(1977), 239-272.

\bibitem {KM}K. Kodaira and Morrow, \textquotedblright\textit{Complex
Manifolds\textquotedblright.}

\bibitem {L}R. Lazarfeld, "\textit{Positivity in Algebraic geometry I and II",
}Springer, 2004.

\bibitem {LSTY}K. Liu, X. Sun, A. Todorov, Shing-Tung Yau, "\textit{The Moduli
Space of} \textit{CY Manifolds} II (Existence of K\"{a}hler-Einstein Metrics
on the Moduli Space of CY Manifolds)", preprint.

\bibitem {LTYZ}K. Liu, A. Todorov, Shing-Tung Yau and K. Zuo,
\textit{"Shafarevich Conjecture for CY Manifolds I" (Moduli of CY Manifolds),
}Quaterly Journal of Pure and Applied Mathematics, vol. \textbf{1.}

\bibitem {ls2}Z. Lu and X. Sun, "T\textit{he Weil-Petersson Volume of the
Moduli Space of Calabi-Yau Manifolds"}. Comm. Math. Phys., \textbf{vol. 261,
}\ No 2, p. 297--322, 2006.

\bibitem {Roe}John Roe, \textquotedblright\textit{Elliptic Operators, Topology
and Asymptotic Methods\textquotedblright} Pitman Research Notes in Mathematics
Series \textbf{179,} Longman Scientific \& Technical, \ 1988.

\bibitem {RS}D. Ray and I. Singer, \textquotedblright\textit{Analytic Torsion
for Complex Manifolds\textquotedblright,} Ann. Math. \textbf{98 }(1973) 154-177.

\bibitem {Sul}D. Sullivan, \textquotedblright\textit{Infinitesimal
Computations in Topology\textquotedblright, Publ. Math. IHES, No \textbf{47
(}1977), 269-331.}

\bibitem {Ti}G. Tian, ''\textit{Smoothness of the Universal Deformation Space
of Calabi-Yau Manifolds and its Petersson-Weil Metric'',} Math. Aspects of
String Theory, ed. S.-T.Yau, World Scientific (1998), 629-346.

\bibitem {To89}A. Todorov, \textit{\textquotedblright The Weil-Petersson
Geometry of \ Moduli Spaces of SU(n}$\geq3$\textit{) (Calabi-Yau Manifolds)
I\textquotedblright,} Comm. Math. Phys. \textbf{126} (1989), 325-346.

\bibitem {To04}A. Todorov, \textit{\textquotedblright Weil-Petersson Volumes
of the Moduli Spaces of CY\ Manifolds", }Communication in Analysis and
geometry, \textbf{Vol. 15, }No2, p. 407-434.

\bibitem {W}E. Viehweg, \textit{\textquotedblright Quasi-Projective Moduli for
Polarized Manifolds\textquotedblright, }\ Ergebnisse der Mathematik und iher
Grenzgebiete 3. Folge, Band \textbf{30}, Springer-Ver;ag, 1991.

\bibitem {Yau}S. T. Yau, \textquotedblright\textit{On the Ricci Curvature of
Compact K\"{a}hler Manifolds and Complex Monge-Amper Equation
I\textquotedblright}, Comm. Pure and App. Math. \textbf{31}(1979)\textbf{,} 339-411.

\bibitem {yau}Shing-Tung Yau, "\textit{A General Schwarz Lemma for Kahler
Manifolds", }American Journal of Mathematics, \textbf{Vol}. 100, No. 1 (Feb.,
1978), pp. 197-203.
\end{thebibliography}
\end{document}